\documentclass[11pt]{article}

\usepackage{graphicx}
\usepackage{latexsym,amsmath,amsfonts,amscd, amsthm, dsfont}
\usepackage{bm,color}
\usepackage{epsfig,verbatim,epstopdf,graphics}
\usepackage{subfigure}
\usepackage{changebar}
\usepackage{multirow}

\usepackage{algorithmic}

\usepackage{yhmath}
 \usepackage{booktabs} 
 \usepackage{tikz}
\usepackage{verbatim}
\usetikzlibrary{arrows,backgrounds,snakes,shapes}
 \numberwithin{equation}{section}

\graphicspath{{./}{./figures/}}
\allowdisplaybreaks

\newtheoremstyle{plainNoItalics}{}{}{\normalfont}{}{\bfseries}{.}{ }{}

\theoremstyle{plain}
\newtheorem{thm}{Theorem}[section]

\theoremstyle{plainNoItalics}

\newtheorem{rem}[thm]{Remark}

\newtheorem{exa}[thm]{Example}

\newcommand{\bE}{{\bf E}}

\newcommand{\bV}{{\bf V}}

\newcommand{\be}{\begin{eqnarray}}
\newcommand{\ee}{\end{eqnarray}}
\newcommand{\beno}{\begin{eqnarray*}}
\newcommand{\eeno}{\end{eqnarray*}}

\makeatletter

\newcommand{\Rmnum}[1]{\expandafter\@slowromancap\romannumeral #1@}
\makeatother

\begin{document}

\baselineskip=1.8pc


\begin{center}
{\bf
Comparison of  semi-Lagrangian discontinuous Galerkin schemes  for
linear and nonlinear transport simulations
}
\end{center}

\vspace{.2in}
\centerline{
Xiaofeng Cai\footnote{
 Department of Mathematical Sciences, University of Delaware, Newark, DE, 19716. E-mail: xfcai@udel.edu.
},
 Wei Guo\footnote{
Department of Mathematics and Statistics, Texas Tech University, Lubbock, TX, 70409. E-mail:
weimath.guo@ttu.edu. Research
is supported by NSF grant NSF-DMS-1620047 and NSF-DMS-1830838.
},
Jing-Mei Qiu\footnote{Department of Mathematical Sciences, University of Delaware, Newark, DE, 19716. E-mail: jingqiu@udel.edu. Research supported by NSF grant NSF-DMS-1522777 and 1818924, Air Force Office of Scientific Computing FA9550-18-1-0257.}
}

\bigskip
\noindent
{\bf Abstract.} Transport problems arise across diverse fields of science and engineering. Semi-Lagrangian (SL) discontinuous Galerkin (DG) methods are a class of high order deterministic transport solvers that enjoy advantages of both SL approach and DG spatial discretization.
In this paper, we review existing SLDG methods to date, and compare numerical their performances. In particular, we make a comparison between the splitting and  non-splitting SLDG methods for multi-dimensional transport simulations. Through extensive numerical results, we offer a practical guide for choosing optimal SLDG solvers for linear and nonlinear transport simulations.

\vfill

\newpage

\section{Introduction}

Semi-Lagrangian (SL) discontinuous Galerkin (DG) methods are a class of high order transport solvers that enjoy the computational advantages of both SL approach and DG spatial discretization. In this paper, we conduct a systematic comparison for  several existing SLDG methods in the literature by considering aspects including accuracy, CPU efficiency, conservation properties, implementation difficulties, with the aim to provide a brief survey of the recent development along this line of research for simulating linear and nonlinear transport problems.

The DG methods belong to the family of finite element methods, which employ piecewise polynomials as approximation and test function spaces. Such methods have undergone rapid development for simulating partial differential equations (PDEs) over the last few decades. For the time-dependent transport simulations, the DG methods are often coupled with the method-of-lines Eulerian framework, by using appropriate time integrators for time evolution, e.g. the well-known TVD Runge-Kutta (RK) methods. It is well-known that the DG method, when coupled with an explicit time integrator suffers very stringent $CFL$ time step restriction for stability, which may be much smaller than that is needed to resolve interesting time scales in physics. Implicit methods can be used to avoid the $CFL$ time step restriction, yet the additional computational cost is required for solving the resulting linear or nonlinear system. The SL approach allows for large time step evolution by incorporating the characteristic tracing mechanism without implicit treatments. Meanwhile, the high order accuracy can be conveniently attained because of its meshed-based nature.  Such a distinguished property makes the SL approach very competitive in transport simulations, and it is becoming more and more popular across diverse fields of science and engineering, such as fluid dynamics \cite{xiu2001semi}, numerical weather prediction \cite{lin1997explicit,Guo2013discontinuous}, plasma physics \cite{cheng1976integration,sonnendrucker1999semi,crouseilles2009conservative}, among many others.
To alleviate the $CFL$ time step restriction associated with the explicit Eulerian DG methods,  the SLDG methods were developed in various settings \cite{restelli2006semi,rossmanith2011positivity,qiu2011positivity,Guo2013discontinuous}. This kind of methods enjoy extraordinary computational advantages for transport simulations, such as high order accuracy, provable unconditional stability and convergence, local mass conservation, small numerical dissipation, superior ability to resolve complex solution structures. Note that the SL framework is highly flexible and hence able to accommodate virtually all the existing spatial discretization methods, see the review paper \cite{qiu2016high}.

In this paper, we are mainly concerned with the SLDG methods. For one-dimensional (1D) transport problems, there exist two types of SLDG methods: the characteristic Galerkin weak formulation and the flux difference formulation. When generalizing to the multi-dimensional case,
similar to other SL methods, the available SLDG methods can be classified into two categories depending on whether the operator splitting strategy is used.
The splitting based SL methods are popular in practice  since one can directly make use of a preferred 1D SL formulation. A celebrated example is the Strang splitting SL method for the Vlasov-Poisson (VP) system proposed by Cheng and Knorr \cite{cheng1976integration}.
In particular, the nonlinearly coupled high-dimensional transport equation is decomposed into a set of  lower dimensional sub-equations that are linear and hence  much easier to evolve numerically in the SL setting. Following a similar idea, many splitting SL schemes based on various spatial discretization are designed, such as
finite volume based \cite{FilbetSB} and finite difference based \cite{Qiu_Christlieb,qiu_shu_sl,qiu2011conservative,cai2016conservative} methods. More recently,  several splitting SLDG methods are proposed in \cite{qiu2011positivity,rossmanith2011positivity}.
On the other hand,  the splitting methods suffer from the inherent splitting error. For instance, if the Strang splitting is used, then the second order splitting error in time is incurred. As mentioned in  \cite{christlieb2014high}, such a relatively low order error may become significant for long term transport simulations and hence greatly compromise the performance of the  SL methods. In addition, for some nonlinear problems, it can be very difficult to track characteristics accurately for the split sub-equations, posing  challenges for the development of high order splitting SL methods. This issue will be discussed in detail in Section \ref{sec:splitting}.
To completely avoid the splitting error, several non-splitting SL methods are developed. Here, we mention the conservative SL multi-tracer transport methods (CSLAM) \cite{lauritzen2010conservative,harris2011flux} and the SL spectral element transport methods \cite{erath2014conservative,Bochev2015conservative}. The first non-splitting SLDG scheme is proposed in \cite{restelli2006semi}. Such a method is based on a
flux difference form but subject to a $CFL$ time step restriction, which degrades its computational efficiency to some extent. Recently, a line of research has been carried out for the development of non-splitting SLDG methods \cite{lee2016high, cai2016high,cai2018high,cai2018Euler}. The proposed methods are unconditionally stable, leading to immense computational efficiency. Meanwhile, the implementation is much more involved than the splitting counterpart.

In this paper, we review the recent development of the state-of-the-art SLDG methods and systematically compare these methods in various settings.
Our goal is to provide a practical reference guide and present extensive numerical results with an emphasis on the study of the accuracy and error versus CPU time. Based on our investigation, we observe that the splitting SLDG method may be preferred due to its simplicity of implementation, especially when small $CFLs$ are needed for accuracy. On the other hand, when large $CFLs$ are used and the splitting error becomes pronounced, the non-splitting SL DG methods are preferred for long time simulations. For example, for the nonlinear incompressible flow simulation, the non-splitting method is more efficient and effective, even though its implementation demands a large amount of human effort. Our simulation results demonstrate that the non-splitting method performs better.

The rest of the paper is organized as follows. In Section 2, we review the existing SLDG methods for solving linear and nonlinear problems. In Section 3, we compare performances of SLDG methods via extensive numerical experiments with an emphasis on efficiency comparison between the splitting and non-splitting formulations. The concluding remark is provided in Section 4.

\section{A Zoo of Semi-Lagrangian discontinuous Galerkin methods}

In this section, we review several SLDG formulations to date. We start with 1D formulations for transport. When extending to multi-dimensions, we discuss both splitting and non-splitting strategies. We will highlight the key differences in the scheme formulations together with implementations, and compare their performances in the numerical section.

\subsection{One-dimensional SLDG methods} \label{review_1d}

Consider the following 1D linear transport equation
\begin{equation}
\label{eq:trans}
u_t
+
(a(x,t)u)_x =0,\quad  x\in\Omega
\end{equation}
with proper initial and boundary conditions. For simplicity, we  consider periodic boundary conditions in this paper.
Here $a(x,t)$ is a velocity field, that could be space and time dependent.
The domain is partitioned by non-overlapping intervals, i.e., $\Omega=\cup_{j}I_j$ where $I_j=[x_{j-\frac12}, x_{j+\frac12} ]$. We denote $\Delta x_j=x_{j+\frac12}-x_{j-\frac12}$ as the length of an element and denote $h=\max_j\Delta x_j$ as the mesh size. We further denote $t^n$ as the discrete $n$-th time level with $\Delta t^n =t^{n+1}-t^n$ being the $n$-th time stepping size.
We define the finite dimensional DG approximation space, $V_h^k = \{ v_h:  v_h|_{I_j} \in P^k(I_j) \}$, where $P^k(I_j)$ denotes the set of polynomials of degree at most $k$ over $I_j$. Next we review several SLDG formulations proposed in the past two decades \cite{restelli2006semi, James, Guo2013discontinuous, cai2016high}.

\subsubsection{SLDG method based on characteristic Galerkin weak formulation}
\label{sec:weak}

Based on the characteristic Galerkin weak formulation \cite{Guo2013discontinuous}, we let the test function $\psi(x,t)$ solve the adjoint problem satisfying the final value problem $\psi(t=t^{n+1}) = \Psi(x)$ with $\forall \Psi \in P^k(I_j)$,
\begin{equation}
\psi_t + a(x,t) \psi_x =0, \ t\in[t^n,t^{n+1}].
\label{final-value}
\end{equation}
Equation \eqref{final-value} is in the convective form, hence the solution stays constant along a characteristic trajectory. Such a strategy has been used in the   ELLAM settings (see, e.g., \cite{celia1990eulerian,herrera1993eulerian,ewing1993eulerian,wang1999ellam}) with continuous finite elements.
It can be shown that \cite{Guo2013discontinuous}
\begin{equation}
\frac{d}{dt}\int_{\widetilde{I}_j(t)} u(x,t)\psi(x,t) dx =0,
\label{property}
\end{equation}
where $\widetilde{I}_j(t)$ is a dynamic interval bounded by characteristics emanating from cell boundaries of $I_j$ at $t=t^{n+1}$.
 \eqref{property} leads to the following SLDG weak formulation. Given $u^n\in V_h^k$, we seek $u^{n+1}\in V_h^k$, such that for $\forall \Psi \in P^k(I_j)$, $j=1,\ldots,m,$
\begin{equation}
\int_{I_j} u^{n+1} \Psi dx = \int_{I_{j}^\star } u(x,t^n) \psi(x,t^n)dx,\quad
\label{integral1}
\end{equation}
where $I_{j}^\star  = [x_{j-\frac12}^\star , x_{j+\frac12}^\star]$ with $x_{j\pm\frac12}^\star$ being the feet of trajectory emanating from $(x_{j\pm\frac12}^\star,t^{n+1})$ at the time level $t^n$.
To update the numerical solution $u^{n+1}$, we need to evaluate the integral on the right-hand side (RHS) of equation \eqref{integral1}. Below we review the procedure proposed in \cite{cai2016high}.
\begin{description}
	\item[(1)] Choose $k+1$ interpolation points $x_{j,q}$, $q=0,\ldots,k$, such as the Gauss-Lobatto points over the interval $I_j$, and locate the feet $x_{j,q}^\star$ by numerically solving  the following trajectory equation:
	\begin{equation}
	\label{eq:final}
	\frac{d x(t) }{dt}  = a( x(t),t ),
	\end{equation}
	with a final value $ x(t^{n+1}) = x_{j,q}$ by a high order numerical integrator. See Figure \ref{schematic_1d_1}.
	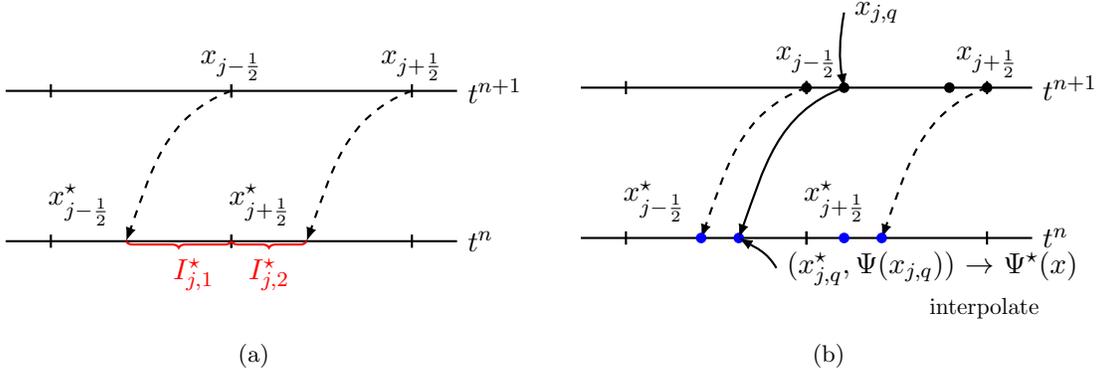
\begin{figure}[h!]
		\centering
		\subfigure[]{
			\begin{tikzpicture}[x=1cm,y=1cm]
			\begin{scope}[thick]
			\draw (-3,3) node[fill=white] {};
			\draw (-3,-1) node[fill=white] {};
			\draw[black]                   (-3,0) node[left] {} -- (3,0)
			node[right]{$t^{n}$};
			\draw[black] (-3,2) node[left] {$$} -- (3,2)
			node[right]{$t^{n+1}$};
			
			\draw[snake=ticks,segment length=2.4cm] (-2.4,2) -- (0,2) node[above] {$x_{j-\frac12}$};
			\draw[snake=ticks,segment length=2.4cm] (0,2) -- (2.4,2) node[above] {$x_{j+\frac12}$};
			
			\draw[snake=ticks,segment length=2.4cm] (-2.4,0) -- (0,0);
			\draw[snake=ticks,segment length=2.4cm] (0,0) -- (2.4,0);
			
			\draw[-latex,dashed](0,2)node[left,scale=1.3]{$$}
			to[out=200,in=70] (-1.4,0) node[above left=2pt] {$x_{j-\frac12}^\star$};  
			\draw[-latex,dashed](2.4,2)node[left,scale=1.3]{$$}
			to[out=200,in=70] ( 1.,0) node[above left=2pt] {$x_{j+\frac12}^\star$};
			
			\draw[snake=brace,mirror snake,red,thick] (-1.4,0) -- (0,0) node[below left=2.5pt] {$I_{j,1}^\star$};
			\draw[snake=brace,mirror snake,red,thick] (0,0) -- (1,0) node[below left=2.5pt] {$I_{j,2}^\star$};
			
			\end{scope}
			\end{tikzpicture}
		}
		\subfigure[]{
			\begin{tikzpicture}[x=1cm,y=1cm]
			\begin{scope}[thick]
			\draw[black]                   (0+1,0) node[left] {} -- (6+1,0)
			node[right]{$t^{n}$};
			\draw[black] (0+1,2) node[left] {$$} -- (6+1,2)
			node[right]{$t^{n+1}$};
			
			\draw[snake=ticks,segment length=2.4cm] (0.6+1,2) -- (3+1,2) node[above] {$x_{j-\frac12}$};
			\draw[snake=ticks,segment length=2.4cm] (3+1,2) -- (5.4+1,2) node[above] {$x_{j+\frac12}$};
			
			\draw[snake=ticks,segment length=2.4cm] (0.6+1,0) -- (3+1,0);
			\draw[snake=ticks,segment length=2.4cm] (3+1,0) -- (5.4+1,0);
			
			\draw[-latex,dashed](3+1,2)node[left,scale=1.3]{$$}
			to[out=200,in=70] (1.6+1,0) node[above left=2pt] {$x_{j-\frac12}^\star$};  
			\draw[-latex,dashed](5.4+1,2)node[left,scale=1.3]{$$}
			to[out=200,in=70] ( 4+1,0) node[above left=2pt] {$x_{j+\frac12}^\star$};

			\fill [black] (3+1,2) circle (2pt) node[] {};
			\fill [black] (5.4+1,2) circle (2pt) node[] {};
			\fill [black] (3+1+0.5,2) circle (2pt) node[] { };
			\fill [black] (5.4+1-0.5,2) circle (2pt) node[] {};
			
			\fill [blue] (3+1 -1.4,0) circle (2pt) node[] {};
			\fill [blue] (5.4+1-1.4,0) circle (2pt) node[] {};
			\fill [blue] (3+1+0.5-1.4,0) circle (2pt) node[] { };
			\fill [blue] (5.4+1-0.5-1.4,0) circle (2pt) node[] {};
			
			\draw[-latex,thick](3+1+0.5,3)node[right,scale=1.0]{$x_{j,q}$}
			to[out=260,in=100] (3+1+0.5,2) node[] {};
			
			\draw[-latex,thick](3+1+0.5,2)node[right,scale=1.0]{}
			to[out=200,in=70] (3+1+0.5-1.4,0) node[] {};
			
			\draw[-latex,thick](3+1+1-1.4,-0.4)node[right=0pt,scale=1.0]{$(x_{j,q}^\star,\Psi(x_{j,q}) )$ $\rightarrow$ $\Psi^\star(x)$ }
			node[below right=8pt,scale=0.8]{$\ \ \ \ \ \ \ \ \ \ \ \ \ \ \ $ interpolate}
			to[out=120,in=330] (3+1+0.5-1.4,0) node[] {};
			
			\end{scope}
			\end{tikzpicture}
		}
	
		\caption{Schematic illustration for the 1D SLDG scheme in \cite{cai2016high}.}	\label{schematic_1d_1}
	\end{figure}

	\item[(2)] Recall that the test function $\psi$ solves the final-value problem \eqref{final-value} and hence stays constant along the characteristics, i.e., $\psi(x_{j,q}^\star ,t^n ) = \Psi(x_{j,q})$.
	Now we are able to determine the unique polynomial $\Psi^\star(x)$ of degree $k$ which interpolates $\psi(x,t^n)$ with $(x_{j,q}^\star , \Psi(x_{j,q}))$ for $q=0,\cdots,k$.

	\item[(3)] Detect intervals/sub-intervals within $I_j^\star = \bigcup_l I_{j,l}^\star$, which are  the intersections between $I_j^\star$ and the grid elements ($l$ is the index for sub-interval). For instance, in Figure \ref{schematic_1d_1} (a), there are two sub-intervals: $I_{j,1}^\star = [x_{j-\frac12}^\star , x_{j-\frac12}] $
	and $I_{j,2}^\star = [x_{j-\frac12},x_{j+\frac12}^\star ]$.

	\item[(4)] Lastly, the RHS of \eqref{integral1} is approximated by
	\begin{equation}
	\int_{I_{j}^\star } u(x,t^n) \psi(x,t^n)dx \approx
	\sum_l \int_{I_{j,l}^\star } u^n(x) \Psi^\star(x)dx.
	\label{integral2}
	\end{equation}
	Note that the integrands in \eqref{integral2} are polynomials of degree $2k$; thus the integration can be evaluated exactly.

\end{description}

%
%
%
%
%
%
%
%

\begin{rem}
There is another slightly different procedure in evaluating the RHS integrals on equation~\eqref{integral1}, that is to locate quadrature points on subintervals, track their values along characteristics curves using equation~\eqref{final-value}, and perform quadrature integrations subinterval by subinterval. For example, see the procedures proposed in \cite{Guo2013discontinuous, lee2016high}. In the case of having constant advection coefficient, the procedure  in \cite{Guo2013discontinuous, lee2016high} and the one reviewed above \cite{cai2016high} are equivalent; they are also equivalent to the shift-projection strategy in the SLDG scheme as proposed in \cite{rossmanith2011positivity}. In the case of the variable coefficient, the procedure reviewed above as proposed in \cite{cai2016high} can be better generalized to the two-dimensional case in Section~\ref{sec: nonsplitting}.
\end{rem}

\subsubsection{SLDG methods based on a flux difference form}
\label{sec: flux_diff}

\begin{figure}[h!]
\centering
\subfigure[]{
\begin{tikzpicture} [scale=1.0]
\centering
  \begin{scope}[thick]
   \draw (-3,3) node[fill=white] {};
    \draw (-3,-1) node[fill=white] {};
    \draw[black]                   (-3,0) node[left] {} -- (3,0)
                                        node[right]{$t^{n}$};
    \draw[black] (-3,3) node[left] {$$} -- (3,3)
                                        node[right]{$t^{n+1}$};

     \draw[snake=ticks,segment length=2.4cm] (-2.4,3) -- (0,3) node[above] {\scriptsize $x_{j-\frac12}$};
     \draw[snake=ticks,segment length=2.4cm] (0,3) -- (2.4,3) node[above] {\scriptsize $x_{j+\frac12}$};

          \draw[snake=ticks,segment length=2.4cm] (-2.4,0) -- (0,0);
     \draw[snake=ticks,segment length=2.4cm] (0,0) -- (2.4,0);

     \draw[-latex,blue](0,0)node[below,scale=1.]{\scriptsize \textcolor{black}{$x_{j-\frac12}$} }
        to[out=90,in=270] (0,3) node[below left=-3pt] { };  

          \draw[-latex,blue](2.4,0)node[below,scale=1.]{\scriptsize \textcolor{black}{$x_{j+\frac12}$} }
        to[out=90,in=270] (2.4,3) node[below left=-3pt] { };

        \fill [red] ( 2.4, 1.5-0.86*1.5+0.1 ) circle (2pt) node[right] {};

         \fill [red] ( 2.4,1.5-0.34*1.5) circle (2pt) node[right] {\scriptsize $(x_{j+\frac12} ,t^{n,s} )$};
          \fill [red] ( 2.4,1.5+0.34*1.5) circle (2pt) node[] {};

           \fill [red] ( 2.4,1.5+0.86*1.5-0.1) circle (2pt) node[] {};

        \fill [red] ( 1.8,0) circle (1.5pt) node[] {};
             \draw[-latex,dashed]( 2.4, 1.5-0.86*1.5+0.1 )node[left,scale=1.3]{$$}
        to[out=200,in=40] ( 1.8,0) node[below=2pt] { };

        \fill [red] ( 0.8,0) circle (1.5pt) node[] {};
        \draw[-latex,dashed]( 2.4, 1.5-0.34*1.5 )node[left,scale=1.3]{$$}
        to[out=200,in=45] ( 0.8,0) node[below=2pt] { };

      \fill [red] (0-0.5,0) circle (1.5pt) node[] {};
        \draw[-latex,dashed]( 2.4, 1.5+0.34*1.5 )node[left,scale=1.3]{$$}
        to[out=195,in=50] ( 0-0.5,0) node[below=2pt] { };

        \fill [red] (-1.5,0) circle (1.5pt) node[] {};
              \draw[-latex,dashed]( 2.4, 1.5+0.86*1.5-0.1 )node[left,scale=1.3]{$$}
        to[out=190,in=55] ( -1.5,0) node[below=2pt] { };

\draw [decorate,color=black,decoration={brace,mirror,amplitude=6pt},xshift=0pt,yshift=0pt]
(0,0) -- (2.4,0) node [black,midway,xshift=0cm,yshift=-12pt]
{\footnotesize $I_{j}$};

        \draw[-latex,thick](0.2,-0.9)node[below right=-8pt,scale=1.0]{\scriptsize $x_{j+\frac12,s}^\star$}
        to[out=80,in=220] (0.8,0) node[] {};

  \end{scope}

\end{tikzpicture}
}
\subfigure[]{
\begin{tikzpicture} [scale=1.0]
\centering
  \begin{scope}[thick]

    \draw[white,fill=blue!3] (2.4,3) to[out=190,in=60] (-1.5,0)  -- (2.4,0)
      -- cycle;
   \node[blue!60, rotate=0] (a) at ( 1. ,1.5) {\LARGE $\Omega_{j+\frac12}$ };
   \draw (-3,3) node[fill=white] {};
    \draw (-3,-1) node[fill=white] {};
    \draw[black]                   (-3,0) node[left] {} -- (3,0)
                                        node[right]{$t^{n}$};
    \draw[black] (-3,3) node[left] {$$} -- (3,3)
                                        node[right]{$t^{n+1}$};

     \draw[snake=ticks,segment length=2.4cm] (-2.4,3) -- (0,3) node[above] {\scriptsize $x_{j-\frac12}$};
     \draw[snake=ticks,segment length=2.4cm] (0,3) -- (2.4,3) node[above] {\scriptsize $x_{j+\frac12}$};

          \draw[snake=ticks,segment length=2.4cm] (-2.4,0) -- (0,0);
     \draw[snake=ticks,segment length=2.4cm] (0,0) -- (2.4,0);

     \draw[-latex,blue](0,0)node[below,scale=1.]{\scriptsize \textcolor{black}{$x_{j-\frac12}$} }
        to[out=90,in=270] (0,3) node[below left=-3pt] { };  

          \draw[-latex,blue](2.4,0)node[below,scale=1.]{\scriptsize \textcolor{black}{$x_{j+\frac12}$} }
        to[out=90,in=270] (2.4,3) node[below left=-3pt] { };

%

             \draw[-latex,blue!60,dashed]( 2.4, 3 )node[left,scale=1.3]{$$}
        to[out=190,in=60] ( -1.5,0) node[below=2pt] { };

       \fill [red] ( 2.4,0) circle (2pt) node[] {};

       \fill [red] ( 0,0) circle (2pt) node[] {};

        \fill [red] ( 1.2+0.44*1.2,0) circle (2pt) node[] {};

        \fill [red] ( 1.2-0.44*1.2,0) circle (2pt) node[] {};

        \fill [red] ( -1.5,0) circle (2pt) node[] {};

        \fill [red] ( -0.75+0.44*0.75,0) circle (2pt) node[] {};

        \fill [red] ( -0.75-0.44*0.75,0) circle (2pt) node[] {};
  \end{scope}

\end{tikzpicture}
}
		
		\caption{Schematic illustration for the 1D SLDG scheme in \cite{restelli2006semi} (a) and the scheme in \cite{qiu2011conservative} (b).}\label{schematic_1d_sl_flux11}
\end{figure}
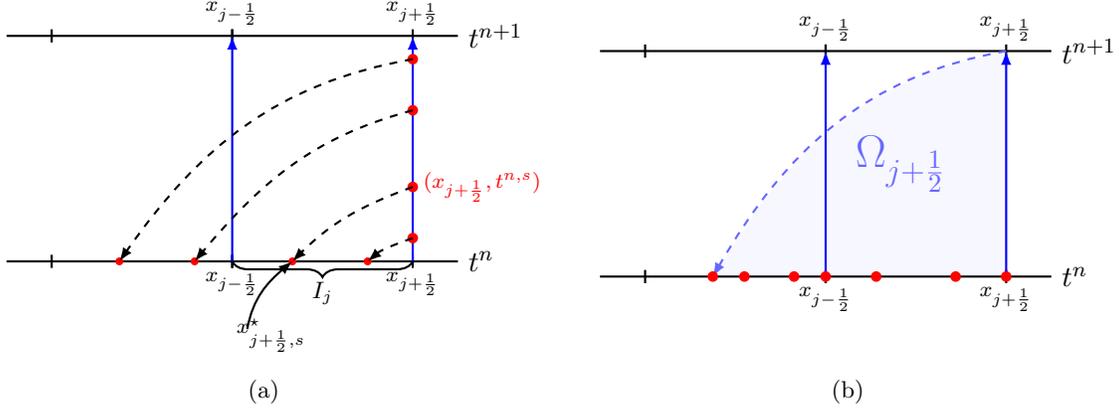

There exist two SLDG methods based on a flux difference form in the literature \cite{restelli2006semi,qiu2011positivity}. Such methods are motivated by the standard DG method for solving hyperbolic conservation laws. Multiplying \eqref{eq:trans} by a time-independent test function $\Psi$, integrating the resulting equation over $I_j$ and performing integration by parts yields
\begin{equation}
\label{eq:semi_dg}
\frac{d}{dt}\int_{I_j}u(x,t)\Psi dx - \int_{I_j}a(x,t)\Psi_xdx + a(x,t)u(x,t)\Psi\Big|_{x_{j+\frac12}}-a(x,t)u(x,t)\Psi\Big|_{x_{j-\frac12}} = 0.
\end{equation}
Instead of coupling a time integrator in the method-of-lines fashion, the SLDG methods make use of the characteristic tracing mechanism. In particular, we further integrate  \eqref{eq:semi_dg} over time interval $[t^n,t^{n+1}]$ and define the SLDG method accordingly. Given data $u^n\in V^k_h$, we seek $u^{n+1}\in V^k_h$ such that for  $\forall\Psi\in V^k_h$
\begin{align}\notag
\int_{I_j}u^{n+1}\Psi\, dx& =\int_{I_j}u^{n}\Psi\, dx  + \int_{t^n}^{t^{n+1}}\left(\int_{I_j}a(x,t)S_t(u^n)\Psi_x\,dx\right)\,dt \\
&+ \int_{t^n}^{t^{n+1}}\left(a(x,t)S_t(u^n)\Psi\Big|_{x_{j+\frac12}^-}\right)\,dt-\int_{t^n}^{t^{n+1}}\left(a(x,t)S_t(u^n)\Psi\Big|_{x_{j-\frac12}^+}\right)\,dt, \label{eq:sldg_flux}
\end{align}
$ n=0,1,\ldots$, where $S_t(u^n)$ denotes the function evolved from the given initial data $u^n$ by following the characteristics, see \cite{restelli2006semi,qiu2011positivity} for more details. The integrals in space and time appearing in \eqref{eq:sldg_flux} are computed by numerical quadrature rules for the SLDG method proposed in \cite{restelli2006semi} on the fully discrete level. For example, to compute the third integral on the RHS of \eqref{eq:sldg_flux}, we locate $p+1$ quadrature nodes   over $[t^n,t^{n+1}]$, denoted as $t^{n,s}$, $s=0,\ldots,p$. Then \[\int_{t^n}^{t^{n+1}}\left(a(x,t)S_t(u^n)\Psi\Big|_{x_{j+\frac12}^-}\right)\,dt\approx\sum_s a(x_{j+\frac12},t^{n,s})S_t(u^n(x^\star_{j+\frac12,s}))\Psi\Big|_{x_{j+\frac12}^-}w_{n,s},
\]
where $x^\star_{j+\frac12,s}$ denote the feet of the characteristic trajectories emanating from the points $(x_{j+\frac12},t^{n,s})$, $s=0,\ldots,p$, see the left panel in Figure \ref{schematic_1d_sl_flux11}, and $w_{n,s}$ denote the associated quadrature weights. As pointed out in \cite{qiu2011positivity}, $S_t(u^n)$ may be a discontinuous function on $[t^n,t^{n+1}]$ for large time step evolution, e.g, $CFL\ge$1, due to the discontinuous nature of $u^n$. In such a case, a  numerical quadrature rule in time may suffer from some stability issue, making the SLDG formulation only conditionally stable. To remedy the drawback, an alternative SLDG method is developed in \cite{qiu2011positivity}, which makes use of the divergence theorem in order to convert the integrals in time into the integrals in space. We still consider the third integral on the RHS of \eqref{eq:sldg_flux} as an example. A direct application of the divergence theorem to the integral form of DG formulation \eqref{eq:semi_dg} over domain $\Omega_{j+\frac12}$ yields
\[
\int_{t^n}^{t^{n+1}}\left(a(x,t)S_t(u^n)\Psi\Big|_{x_{j+\frac12}^-}\right)\,dt = \int_{x^\star_{j+\frac12}}^{x_{j+\frac12}}u^n(x)\,dx\cdot\Psi\Big|_{x_{j+\frac12}^-}.
\]
Hence, numerical quadrature rules are applied to the integrals in space only and in the subinterval-by-subinterval manner, see the right panel in Figure \ref{schematic_1d_sl_flux11}, leading to an unconditionally stable SLDG method.

\begin{rem} Since the SLDG method in \cite{restelli2006semi} is only conditionally stable, the efficiency is limited to some extent. On the other hand, such a formulation can be directly generalized to the multi-dimensional case without operator splitting. For the SLDG method with the divergence theorem \cite{qiu2011conservative}, its direct generalization to a multi-dimensional problem is much more difficult and technically involved, thus dimensional splitting is often used as discussed in Section~\ref{sec:splitting}.
\end{rem}
\begin{rem} It can be shown that the SLDG formulation \eqref{integral1} is equivalent to \eqref{eq:sldg_flux} if all the integrals are computed exactly. In general, we have to resort to numerical quadrature rules or polynomial interpolation, and hence these SLDG formulations have slightly different performance in simulations. 	
\end{rem}

\subsection{Two-dimensional SLDG methods with operator splitting} \label{sec:splitting}

Next we discuss the extension of the 1D SLDG formulations to the two-dimensional (2D) case via dimensional splitting. Consider the following 2D transport equation
\begin{equation}
u_t +
(a(x,y,t)u )_x+
(b(x,y,t)u )_y= 0,
\ (x,y)\in  \Omega,
\label{2D_transport}
\end{equation}
with a proper initial condition $u(x,y,0) = u_0(x,y)$.
Here, $(a(x,y,t),b(x,y,t) )$ is a prescribed velocity field depending on time and space.
Boundary conditions are periodic for simplicity. The rectangular domain $\Omega$ is partitioned in terms of  a Cartesian mesh with each computational cell denoted as $A_{ij} = [x_{i-\frac12}, x_{i+\frac12}]\times[y_{j-\frac12}, y_{j+\frac12}]$. In the dimensional splitting setting, we use the piecewise $Q^k$ tensor-product polynomial spaces. Hence, there are $(k+1)^2$ degrees of freedom per computational cell.

Note that all the 1D SLDG formulations introduced in section \ref{review_1d} can be extended to the 2D case via operator splitting.
Below, we present a general framework of  splitting SLDG algorithms.

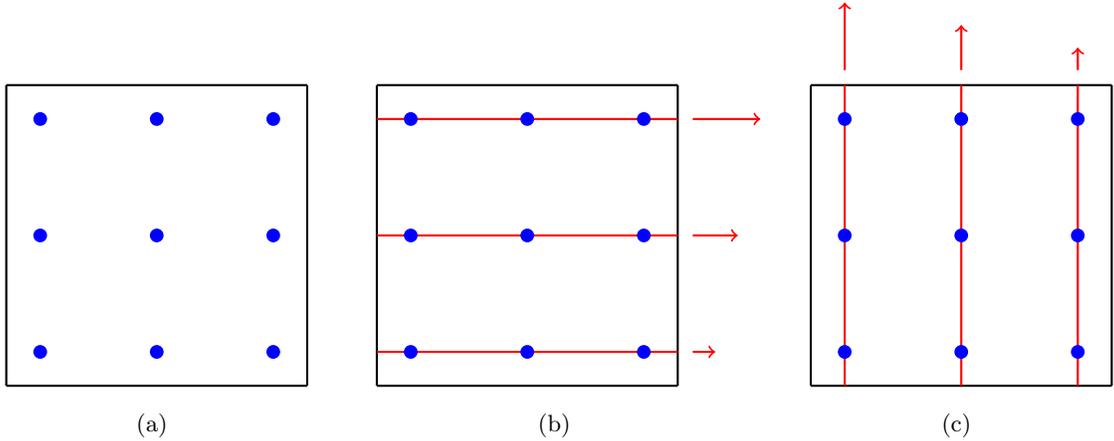
\begin{figure}[h!]
\centering
\subfigure[]{
\begin{tikzpicture}[x=1cm,y=1cm]
    \draw[black,thick] (0,0) node[left] {} -- (4,0)
                                        node[right]{};
    \draw[black,thick] (4,0) node[left] {$$} -- (4,4)
                                        node[right]{};
    \draw[black,thick] (4,4) node[left] {$$} -- (0,4)
                                        node[right]{};
    \draw[black,thick] (0,4 ) node[left] {$$} -- (0,0)
                                        node[right]{};
%
%
     \fill [blue] (2-0.7746*2,2- 0.7746*2 ) circle (2.6pt) node[] {};
     \fill [blue] (2-0.*2,2- 0.7746*2 ) circle (2.6pt) node[] {};
    \fill [blue] (2+0.7746*2,2- 0.7746*2 ) circle (2.6pt) node[] {};

     \fill [blue] (2-0.7746*2,2  ) circle (2.6pt) node[] {};
     \fill [blue] (2-0.*2,2  ) circle (2.6pt) node[] {};
    \fill [blue] (2+0.7746*2,2  ) circle (2.6pt) node[] {};

     \fill [blue] (2-0.7746*2,2+0.7746*2 ) circle (2.6pt) node[] {};
     \fill [blue] (2-0.*2,2+ 0.7746*2 ) circle (2.6pt) node[] {};
    \fill [blue] (2+0.7746*2,2+ 0.7746*2 ) circle (2.6pt) node[] {};
%
%
%
%
%
\end{tikzpicture}
}
\subfigure[]{
\begin{tikzpicture}[x=1cm,y=1cm]
    \draw[black,thick] (0,0) node[left] {} -- (4,0)
                                        node[right]{};
    \draw[black,thick] (4,0) node[left] {$$} -- (4,4)
                                        node[right]{};
    \draw[black,thick] (4,4) node[left] {$$} -- (0,4)
                                        node[right]{};
    \draw[black,thick] (0,4 ) node[left] {$$} -- (0,0)
                                        node[right]{};

    \draw[red,thick] (0,2-0.7746*2) node[left] {} -- (4,2-0.7746*2)
                                        node[right]{};
                                        \draw[->,red,thick] (4.2,2-0.7746*2) --(4.5,2-0.7746*2);
\draw[red,thick] (0,2+0.*2) node[left] {} -- (4,2+0. *2)
                                        node[right]{};
                                        \draw[->,red,thick] (4.2,2+0. *2) --(4.8,2+0. *2);

    \draw[red,thick] (0,2+0.7746*2) node[left] {} -- (4,2+0.7746*2)
                                        node[right]{};
                                        \draw[->,red,thick] (4.2,2+0.7746*2) --(5.1,2+0.7746*2);

     \fill [blue] (2-0.7746*2,2- 0.7746*2 ) circle (2.6pt) node[] {};
     \fill [blue] (2-0.*2,2- 0.7746*2 ) circle (2.6pt) node[] {};
    \fill [blue] (2+0.7746*2,2- 0.7746*2 ) circle (2.6pt) node[] {};

     \fill [blue] (2-0.7746*2,2  ) circle (2.6pt) node[] {};
     \fill [blue] (2-0.*2,2  ) circle (2.6pt) node[] {};
    \fill [blue] (2+0.7746*2,2  ) circle (2.6pt) node[] {};

     \fill [blue] (2-0.7746*2,2+0.7746*2 ) circle (2.6pt) node[] {};
     \fill [blue] (2-0.*2,2+ 0.7746*2 ) circle (2.6pt) node[] {};
    \fill [blue] (2+0.7746*2,2+ 0.7746*2 ) circle (2.6pt) node[] {};

%
%

%
\end{tikzpicture}
}
\subfigure[]{
\begin{tikzpicture}[x=1cm,y=1cm]
       \draw[black,thick] (0,0) node[left] {} -- (4,0)
                                        node[right]{};
    \draw[black,thick] (4,0) node[left] {$$} -- (4,4)
                                        node[right]{};
    \draw[black,thick] (4,4) node[left] {$$} -- (0,4)
                                        node[right]{};
    \draw[black,thick] (0,4 ) node[left] {$$} -- (0,0)
                                        node[right]{};

        \draw[red,thick] (2-0.7746*2,4 ) node[left] {$$} -- (2-0.7746*2,0)
                                        node[right]{};
                                        \draw[->,red,thick] (2-0.7746*2,4.2) --(2-0.7746*2,5.1);
        \draw[red,thick] (2-0. *2,4 ) node[left] {$$} -- (2-0. *2,0)
                                        node[right]{};
                                        \draw[->,red,thick] (2-0. *2,4.2) --(2-0. *2,4.8);
        \draw[red,thick] (2+0.7746*2,4 ) node[left] {$$} -- (2+0.7746*2,0)
                                        node[right]{};
                                        \draw[->,red,thick] (2+0.7746*2,4.2) --(2+0.7746*2,4.5);

     \fill [blue] (2-0.7746*2,2- 0.7746*2 ) circle (2.6pt) node[] {};
     \fill [blue] (2-0.*2,2- 0.7746*2 ) circle (2.6pt) node[] {};
    \fill [blue] (2+0.7746*2,2- 0.7746*2 ) circle (2.6pt) node[] {};

     \fill [blue] (2-0.7746*2,2  ) circle (2.6pt) node[] {};
     \fill [blue] (2-0.*2,2  ) circle (2.6pt) node[] {};
    \fill [blue] (2+0.7746*2,2  ) circle (2.6pt) node[] {};

     \fill [blue] (2-0.7746*2,2+0.7746*2 ) circle (2.6pt) node[] {};
     \fill [blue] (2-0.*2,2+ 0.7746*2 ) circle (2.6pt) node[] {};
    \fill [blue] (2+0.7746*2,2+ 0.7746*2 ) circle (2.6pt) node[] {};
%
%
%
%
%
\end{tikzpicture}
}
\caption{Schematic illustration of the 2D SLDG scheme via Strang splitting. $k=2$.}\label{splitting}
\end{figure}

\begin{description}
  \item[(1)] Locate $(k+1)^2$ tensor-product Gaussian nodes on cell $A_{ij}$:  $(x_{i,p},y_{j,q})$, $p,q=0,\cdots,k$.
      For example, see Figure \ref{splitting} (left) for the case of $k=2$.

  \item[(2)]  Split the equation \eqref{2D_transport} into two 1D advection problems based on the quadrature nodes in both $x$- and $y$- directions:
        \begin{align}
         u_t +
         (a(x,y,t)u )_x = 0, \label{split_x} \\
   u_t +
  (b(x,y,t)u )_y = 0. \label{split_y}
        \end{align}
  \item[(3)]  Based on a 1D SLDG formulation, evolve the split equations \eqref{split_x}-\eqref{split_y} via Strang splitting over a time step $\Delta t$: (a) Evolve 1D equation \eqref{split_x} at each quadrature node $y_{j,q}$ for a half time-step $\Delta t/2$, see Figure \ref{splitting} (middle); (b)  Evolve 1D equation \eqref{split_y} at each quadrature node $x_{i,p}$ for a full time-step $\Delta t$, see Figure \ref{splitting} (right); (c) Evolve 1D equation \eqref{split_x} at each quadrature node $y_{j,q}$ for another half time-step $\Delta t/2$, see Figure \ref{splitting} (middle).

\end{description}

\begin{rem}
The splitting 2D SLDG formulation inherits many desired properties from the base 1D formulation, such as high order accuracy in space, unconditional stability, and mass conservation. Meanwhile, a splitting error in time is introduced. For the Strang splitting, the error is of second order. Higher order splitting methods can be constructed in the spirit of composition methods \cite{yoshida, hairerGeom}, while the number of intermediate stages, hence the computational cost, increases exponentially with the order of the splitting method. For example, a fourth order splitting SLDG method is developed in \cite{rossmanith2011positivity} for solving the VP system.
\end{rem}

\subsection{Two-dimensional SLDG method without operator splitting}\label{sec: nonsplitting}

Recently, a class of high order non-splitting SLDG methods is under great development \cite{cai2016high,cai2018high,cai2018Euler} and in \cite{lee2016high} for unstructured meshes. They are unconditionally stable and mass conservative. Below we briefly describe the 2D SLDG algorithm proposed in \cite{cai2016high}, which is a direct generalization of 1D algorithm in Section \ref{review_1d}.
We assume that the domain $\Omega$ is rectangular and is partitioned by a set of non-overlapping tensor-product rectangular elements $A_j$, $j=1,\ldots,J$. We remark that, unlike the splitting formulation, such a non-splitting SLDG method is based on a truly multi-dimensional formulation and hence allows for an unstructured mesh \cite{lee2016high}. We then define the finite dimensional DG approximation space, $\bV_h^k = \{ v_h:  v_h|_{A_j} \in P^k(A_j) \}$, where $P^k(A_j)$ denotes the set of polynomials of degree at most $k$ over $A_j$.


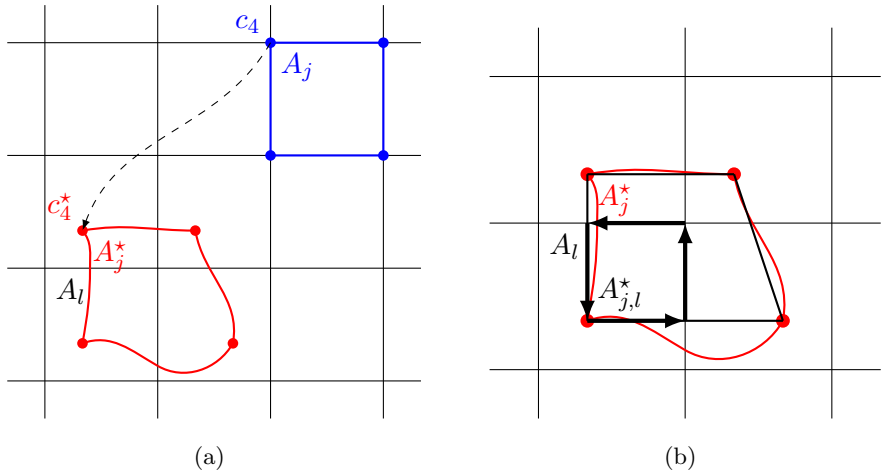
\begin{figure}[h]
\centering
\subfigure[]{
\begin{tikzpicture}
    \draw[black,thin] (0,0.5) node[left] {} -- (5.5,0.5)
                                        node[right]{};
    \draw[black,thin] (0,2.) node[left] {$$} -- (5.5,2)
                                        node[right]{};
    \draw[black,thin] (0,3.5) node[left] {$$} -- (5.5,3.5)
                                        node[right]{};
    \draw[black,thin] (0,5 ) node[left] {$$} -- (5.5,5)
                                        node[right]{};
    \draw[black,thin] (0.5,0) node[left] {} -- (0.5,5.5)
                                        node[right]{};
    \draw[black,thin] (2,0) node[left] {$$} -- (2,5.5)
                                        node[right]{};
    \draw[black,thin] (3.5,0) node[left] {$$} -- (3.5,5.5)
                                        node[right]{};
    \draw[black,thin] (5,0) node[left] {$$} -- (5,5.5)
                                        node[right]{};
    \fill [blue] (3.5,3.5) circle (2pt) node[] {};
    \fill [blue] (5,3.5) circle (2pt) node[] {};
    \fill [blue] (3.5,5) circle (2pt) node[below right] {$A_j$} node[above left] {$c_4$};
    \fill [blue] (5,5) circle (2pt) node[] {};

     \draw[thick,blue] (3.5,3.5) node[left] {} -- (3.5,5)
                                        node[right]{};
      \draw[thick,blue] (3.5,3.5) node[left] {} -- (5,3.5)
                                        node[right]{};
       \draw[thick,blue] (3.5,5) node[left] {} -- (5,5)
                                        node[right]{};
        \draw[thick,blue] (5,3.5) node[left] {} -- (5,5)
                                        node[right]{};
    \fill [red] (1.,1) circle (2pt) node[above right,black] {};
    \fill [red] (3,1) circle (2pt) node[] {};
    \fill [red] (1,2.5) circle (2pt) node[below right] {$A_j^\star$} node[above left] {$c_4^\star$};
    \fill [red] (2.5,2.5) circle (2pt) node[] {};

     \draw[-latex,dashed](3.5,5)node[right,scale=1.0]{}
        to[out=240,in=70] (1,2.50) node[] {};

     \draw (0.5+0.01,2-0.01) node[fill=white,below right] {$A_l$};

     \draw [red,thick] (1,1)node[right,scale=1.0]{}
        to[out=20,in=150] (2,0.7) node[] {};

        \draw [red,thick] (2,0.7)node[right,scale=1.0]{}
        to[out=330,in=240] (3,1) node[] {};
             \draw [red,thick] (1,2.5)node[right,scale=1.0]{}
        to[out=310,in=90] (1.1,2) node[] {};
        \draw [red,thick] (1.1,2)node[right,scale=1.0]{}
        to[out=270,in=80] (1,1) node[] {};

        \draw [red,thick] (1,2.5)node[right,scale=1.0]{}
        to[out=10,in=180] (2.5,2.5) node[] {};

        \draw [red,thick] (3,1)node[right,scale=1.0]{}
        to[out=80,in=280] (2.5,2.5) node[] {};
\end{tikzpicture}
}
\subfigure[]{

\begin{tikzpicture}[scale = 1.3]
    \draw[black,thin] (0,0.5) node[left] {} -- (4,0.5)
                                        node[right]{};
    \draw[black,thin] (0,2.) node[left] {$$} -- (4,2)
                                        node[right]{};
    \draw[black,thin] (0,3.5) node[left] {$$} -- (4,3.5)
                                        node[right]{};
    \draw[black,thin] (0.5,0) node[left] {} -- (0.5,4)
                                        node[right]{};
    \draw[black,thin] (2,0) node[left] {$$} -- (2,4)
                                        node[right]{};
    \draw[black,thin] (3.5,0) node[left] {$$} -- (3.5,4)
                                        node[right]{};

    \fill [red] (1.,1) circle (2pt) node[above right,black] {$A_{j,l}^{\star}$};
    \fill [red] (3,1) circle (2pt) node[] {};
    \fill [red] (1,2.5) circle (2pt) node[below right] {$A_j^{\star}$} node[above left] {};
    \fill [red] (2.5,2.5) circle (2pt) node[] {};

     \draw (0.5+0.01,2-0.01) node[fill=white,below right] {};

\draw [red,thick] (1,1)node[right,scale=1.0]{}
        to[out=20,in=150] (2,0.7) node[] {};

        \draw [red,thick] (2,0.7)node[right,scale=1.0]{}
        to[out=330,in=240] (3,1) node[] {};
             \draw [red,thick] (1,2.5)node[right,scale=1.0]{}
        to[out=310,in=90] (1.1,2) node[] {};
        \draw [red,thick] (1.1,2)node[right,scale=1.0]{}
        to[out=270,in=80] (1,1) node[] {};

        \draw [red,thick] (1,2.5)node[right,scale=1.0]{}
        to[out=10,in=180] (2.5,2.5) node[] {};

        \draw [red,thick] (3,1)node[right,scale=1.0]{}
        to[out=80,in=280] (2.5,2.5) node[] {};
           \draw (0.5+0.01,2-0.01) node[fill=white,below right] {$A_l$};
         \draw[-latex,ultra thick] (1,1)node[right,scale=1.0]{}
        to  (2,1) node[] {};

             \draw[-latex,ultra thick]  (2,1)node[right,scale=1.0]{}
        to (2,2) node[] {};
             \draw[-latex,ultra thick]  (2,2)node[right,scale=1.0]{}
        to (1,2) node[] {};
             \draw[-latex,ultra thick]   (1,2)node[right,scale=1.0]{}
        to (1,1) node[] {};

\draw [thick] (1,1)-- (3,1) node[] {};
 \draw [thick] (1,2.5) -- (1,1) node[] {};
        \draw [thick] (1,2.5)--(2.5,2.5) node[] {};
        \draw [thick] (3,1)--(2.5,2.5) node[] {};
\end{tikzpicture}

}
\caption{Schematic illustration of the SLDG formulation in two dimension: quadrilateral approximation to an upstream cell.  }
\label{schematic_2d}
\end{figure}

Similar to the 1D case, we need the weak formulation of characteristic Galerkin formulation \cite{Guo2013discontinuous,cai2016high}.
Specifically, let the test function $\psi(x,y,t)$ satisfy the adjoint problem with   $\Psi\in P^k(A_j)$,
\begin{equation}
\begin{cases}
\psi_t + a(x,y,t) \psi_x + b(x,y,t) \psi_y = 0,\\
\psi(t=t_{n+1}) = \Psi.
\end{cases}
\end{equation}
Then, we have the identity
\begin{equation}
\frac{d}{dt} \int_{\widetilde{A}_{j}(t) } u(x,y,t) \psi(x,y,t) dxdy =0,
\label{2d_dt}
\end{equation}
where $\widetilde{A}_{j}(t)$ is the dynamic cell, moving from the Eulerian cell $A_{j}$ at $t^{n+1}$, i.e.,  $A_j=\widetilde{A}_{j}(t^{n+1})$, backward in time by following the characteristics trajectories. We denote $\widetilde{A}_{j}(t^n)$ as $A_{j}^\star$, i.e., the upstream cell bounded by the red curves in Figure \ref{schematic_2d}. The SLDG method is defined as follows. Given $u^n\in \bV_h^k$, we seek $u^{n+1}\in \bV_h^k$, such that for $\forall \Psi\in P^k(A_j)$, $j=1,\ldots,J$,
\begin{equation}
\int_{ A_{j}  } u^{n+1} \Psi dxdy =
\int_{ A_j^\star } u^n\psi(x,y,t^{n} ) dxdy.
\label{sl_2d}
\end{equation}
As in the 1D case, the main difficulty lies in the evaluation of the RHS in \eqref{sl_2d}. Below, we briefly discuss the procedure.




Denote $c_q$, $q=1,\cdots,4$ as the four vertices of $A_j$  with the coordinates $(x_{j,q} , y_{j,q} )$.
We trace characteristics backward in time to $t^n$ for the four vertices
 by numerically solving the characteristics equations,
\begin{equation}
\begin{cases}
\frac{d x(t) }{dt} = a(x(t) ,y(t) ,t ),\\
\frac{ d y(t) }{dt} = b(x(t) ,y(t) ,t ), \\
x(t^{n+1} ) = x_{j,q},\\
y(t^{n+1} ) = y_{j,q},
\end{cases}
\label{characteristic}
\end{equation}
 and  obtain $c_q^\star$ with the new coordinate $( x_{j,q}^\star , y_{j,q}^\star ),\,q=1,\cdots,4$.
For example, see $c_4$ and $c_4^\star$ in Figure \ref{schematic_2d}.
The upstream cell $A_j^\star$ can be approximated by a quadrilateral determined by the four vertices $c_q^\star$. Note that $u^n$ is a piecewise polynomial based on the partition. Then the integral over $A_j^\star$ has to be evaluated subregion-by-subregion. To this end, we denote $A_{j,l}^\star$ as a non-empty overlapping region between the upstream cell $A_j^\star$ and the background grid cell $A_l$, i.e., $A_{j,l}^\star = A_{j}^\star \cap A_l$, $A_{j,l}^\star \neq \emptyset$, and define the index set $\varepsilon_j^\star:=\{ l| A_{j,l}^\star \neq \emptyset \}$, see Figure \ref{schematic_2d} (b). The detailed procedure of detecting $A_{j,l}^\star$ can be found in \cite{cai2016high}.
The integral over the upstream cell $A_{j}^\star$ is broken up into the following integrals,
\begin{equation}
\int_{ A_{j}  } u^{n+1} \Psi dxdy
=
\sum_{l\in\varepsilon_j^\star } \int_{ A_{j,l}^\star } u^n\psi(x,y,t^{n} ) dxdy.
\label{temp1}
\end{equation}
Furthermore, $\psi(x,y,t^n)$ is not a polynomial in general, posing additional challenges for evaluating the integrals on the RHS of \eqref{temp1}. On the other hand, if the velocity field $(a,b)$ is smooth, then $\psi(x,y,t^n)$ is smooth accordingly and can be well approximated by a polynomial. The following  procedure is then proposed.

 \begin{enumerate}
   \item  \textit{Least-squares approximation of test function $\psi(x,y,t^n)$.}
We use a least-squares procedure to approximate the test function $\psi(x,y,t^n)$  by a polynomial, based on the fact that  $\psi$ stays constant along characteristics. In particular, for $k=1$, we reconstruct a $P^1$ polynomial $\Psi^\star(x,y)$ by least-squares  with the interpolation constraints
\[
\Psi^\star(x_{j,q}^\star,y_{j,q}^\star) = \Psi(x_{j,q},y_{j,q}),\quad q=1,\ldots,4.
\]
Then,
\begin{equation}
\sum_{l\in\varepsilon_j^\star } \int_{ A_{j,l}^\star } u^n\psi(x,y,t^{n} ) dxdy \approx \sum_{l\in\varepsilon_j^\star } \int_{ A_{j,l}^\star } u^n\Psi^\star(x,y) dxdy.
\label{temp2}
\end{equation}


   \item  \textit{Line integral evaluation via Green's theorem.}
   Note that the integrands on the RHS of \eqref{temp2} are piecewise polynomials. To further simplify the implementation, we make use of  Green's theorem. We first introduce two auxiliary polynomial functions $P(x,y)$ and $Q(x,y)$ such that
   \begin{equation*}
   -\frac{\partial P }{\partial y } + \frac{\partial Q}{\partial x }  =  u(x,y,t^n)\Psi^\star(x,y).
   \end{equation*}
   Then area integral $ \int_{A_{j,l}^\star } u^n\Psi^\star(x,y)dxdy  $ can be converted into line integrals via Green's theorem, i.e.,
   \begin{equation}
   \int_{A_{j,l}^\star } u(x,y,t^n)\Psi^\star(x,y)dxdy = \oint_{\partial A_{j,l}^\star}  Pdx + Qdy,
   \label{Green}
   \end{equation}
   see Figure \ref{schematic_2d} (b). Note that the choices of $P$ and $Q$ are not unique, but the value of the line integrals is independent of the choices. For the implementation details, see \cite{cai2016high}.
 \end{enumerate}

\begin{rem}  The quadrilateral approximation to upstream cells is second order accurate. To achieve third order accuracy, one can use the quadratic-curved quadrilateral approximation, see \cite{cai2016high} for more details.
\end{rem}

\subsection{SLDG schemes for nonlinear models}
\label{sec: splitting}

In this subsection, we discuss the application of the SLDG methods  in nonlinear modelings, e.g. the nonlinear VP system,  the guiding center Vlasov equation, and the incompressible Euler equation.

\begin{description}
  \item[1.] \emph{Vlasov-Poisson system.}
  We consider the 1D1V VP system,
\begin{equation}
f_t + vf_x + E(x,t) f_v =0,
\label{equation:vp}
\end{equation}
\begin{equation}
E(x,t) = -\phi_x, \
-\phi_{xx}(x,t) = \rho(x,t).
\label{eq:poisson}
\end{equation}
Here $x$ and $v$ are the coordinates in the phase space $(x,v)\in \Omega_x \times \mathbb{R}$,
$f(x,v,t)$ is the probability distribution function describing the probability of finding a particle with velocity $v$ at position $x$ and at time $t$. Note that in the Vlasov equation \eqref{equation:vp}, $f$ is accelerated by the electric field $E$ which is determined by Poisson's equation  \eqref{eq:poisson}. $\phi$ is the self-consistent electrostatic potential, and  $\rho(x,t)= \int_{\mathbb{R}} f(x,v,t) dv -1$ denotes charge density.
Here, we assume that infinitely massive ions are uniformly distributed
in the background.
We recall several conserved quantities in the VP system below, which should remain constant in time:
\begin{itemize}
  \item $L^p$ norm, $1\leq p \leq \infty$:
      \begin{equation}
        \| f \|_p = \left(  \int_v \int_x |f(x,v,t)|^p\,dxdv \right)^{\frac{1}{p}},
      \end{equation}
  \item Energy:
      \begin{equation}
        \text{Energy} = \int_v \int_x f(x,v,t) v^2 dx dv + \int_x E^2 (x,t)dx,
      \end{equation}
  \item Entropy:
     \begin{equation}
      \text{Entropy} = \int_v \int_x f(x,v,t) \log (f(x,v,t) ) dxdv.
     \end{equation}
\end{itemize}

  \item[2.] \emph{Incompressible Euler equations and the guiding center Vlasov model.}

The 2D time dependent incompressible Euler equations in the vorticity-stream function formulation reads
\begin{equation}
\begin{split}
\omega_t + \nabla\cdot (\mathbf{u}\omega) = 0, \\
\Delta \Phi = \omega,\
\mathbf{u}
=
( -\Phi_y , \Phi_x),
\end{split}
\label{Euler}
\end{equation}
where $\mathbf{u}$ is the velocity field, $\omega$ is the vorticity of the fluid, and $\Phi$ is the stream-function determined by Poisson's equation. The other closely related model concerned is the guiding center approximation of the 2D Vlasov model, which describes a highly magnetized plasma in the
transverse plane of a tokamak \cite{shoucri1981two,crouseilles2009conservative,frenod2015long,yang2014conservative}:
\begin{align}\label{guiding_center}
\rho_t + \nabla\cdot (\mathbf{E}^\perp\rho) = 0, \\
-\Delta \Phi = \rho,\ \bE^\perp = ( -\Phi_y , \Phi_x)
, \label{poisson}
\end{align}
where the unknown variable $\rho$ denotes the charge density of the plasma, and $\bE$, determined by $\bE = - \nabla \Phi$, is the electric field.
Despite their different application backgrounds, the above two models indeed have an equivalent mathematical formulation up to a sign difference in Poisson's equation. For both models,  the following physical quantities remain constant over time
  \begin{itemize}
   \item  Energy:
            \begin{equation*}
            \| \mathbf{u} \|_{L^2}^2 = \int_{\Omega} \mathbf{u}\cdot\mathbf{u}dxdy, \quad \mbox{(Euler),}
\quad
            \| \mathbf{E} \|_{L^2}^2 = \int_{\Omega} \mathbf{E}\cdot\mathbf{E}dxdy, \quad \mbox{(Vlasov).}
            \end{equation*}
   \item   Enstrophy:
            \begin{equation*}
            \| \omega  \|_{L^2}^2 = \int_{\Omega} \omega^2 dxdy,  \quad \mbox{(Euler),}
\quad            \| \rho  \|_{L^2}^2 = \int_{\Omega} \rho^2 dxdy, \quad \mbox{(Vlasov).}
            \end{equation*}
 \end{itemize}
\end{description}





The splitting SLDG methods for linear transport can be directly applied to both nonlinear models. Below, we summarize the procedures. The notation is consistent with that used in Section \ref{sec:splitting}.
\begin{description}
	\item[(1)]
  For the VP system \eqref{equation:vp}, we follow the classic SL approach with Strang splitting   \cite{cheng1976integration,rossmanith2011positivity,qiu2011conservative,qiu2011positivity,einkemmer2017study}.
 Over one time step $\Delta t$, we perform the following steps:
\begin{enumerate}
  \item Evolve split 1D equation $f_t+vf_x=0$ at each quadrature node $v_{j,q}$ for a half time-step $\Delta t/2$.
  \item Solve Poisson's equation $-\phi_{xx}(x,t^{n+\frac12} )=\rho(x,t^{n+\frac12} )$, and compute $E^{n+\frac12} = -\phi_x(x,t^{n+\frac12} )$.
  \item Evolve split 1D equation $f_t+ E^{n+\frac12} f_v =0$ at each quadrature node $x_{i,p}$ for a full time-step $\Delta t$.
  \item Evolve 1D equation $f_t+vf_x=0$ at each quadrature node $v_{j,q}$ for a half time-step $\Delta t/2$.
\end{enumerate}
The Poisson's equation above can be solved by any proper elliptic solvers, such as the LDG method \cite{arnold2002unified,rossmanith2011positivity}. The Strang splitting strategy is subject to a second order splitting error in time.

\item[(2)] If we directly apply the Strang splitting idea to the 2D incompressible Euler equation in the vorticity-stream function formulation as well as the guiding center Vlasov model,  the procedures are outlined below
\begin{enumerate}
  \item  Obtain $\Phi_y$ by solving $\Delta \Phi  = \omega^n$, where  $\omega^n$ denotes the numerical solution of $\omega(x,y,t^n)$.
  \item  Evolve split 1D equation $\omega_t -( \Phi_y \omega )_x =0$ at each quadrature node  $y_{j,q}$ for a half time-step $\Delta t/2$ and then obtain  solution $\omega^\star$.
  \item Obtain $\Phi_x$ by solving $\Delta \Phi  = \omega^\star$.
  \item Evolve 1D equation $\omega_t + (\Phi_x \omega)_y =0$ at each quadrature node $x_{i,p}$ for a full time-step $\Delta t$ and  obtain  solution $\omega^{\star\star}$.
  \item Obtain $\Phi_y$ by solving $\Delta \Phi  = \omega^{\star\star}$.
  \item Evolve 1D equation $\omega_t -( \Phi_y \omega )_x =0$ at each quadrature node $y_{j,q}$ for a half time-step $\Delta t/2$, and obtain the solution $\omega^{n+1}$ on next time step $t^{n+1}$.
\end{enumerate}


\end{description}
\begin{rem}
Despite the Strang splitting procedure, the above algorithm is only of first order accuracy in time. The reason is that the field $\Phi$ is time-dependent for the split 1D equation. When the field equation and transport equation are solved in a step-by-step manner, numerical tracing of characteristics feet is subject to a first order error in time, which is difficult to enhance beyond first order accuracy. Note that for the case of the VP system, the field equation $\phi$ and hence the electric field $E$ does not change in time, when the \emph{split} 1D equation $f_t + E f_v =0$ is evolved.
\end{rem}

When applying the 2D SLDG methods for the nonlinear Vlasov and incompressible Euler models, the main difficulty is to track nonlinear characteristics with high order temporal accuracy. We adopt the two-stage multi-derivative prediction-correction algorithms as proposed in \cite{qiu2017high} for the VP system and in \cite{xiong2018high} for the Euler equation and the guiding center model. We refer the reader to \cite{cai2018high,cai2018Euler} for details.

\section{Numerical tests} \label{numerical}

In this section, we perform extensive numerical tests to compare the SLDG methods reviewed in the previous section for solving linear and nonlinear transport problems. Note that the 1D SLDG methods based on the characteristic Galerkin formulation in Section~\ref{sec:weak} and based on the flux difference form with unconditional stability in Section~\ref{sec: flux_diff} perform virtually the same in terms of accuracy and CPU cost. Hence, we only apply the SLDG formulation in Section~\ref{sec:weak} with dimensional splitting for its comparison against the non-splitting SLDG method for 2D transport simulations. Examples include 2D linear passive-transport problems in Section \ref{num_transport}, the nonlinear VP system in Section~\ref{num_VP}, and  the incompressible Euler equation and the guiding center Vlasov model in Section \ref{num_Euler}. The performance of the two SLDG methods are benchmarked in terms of {\em error magnitude, order of convergence in space and time, CPU cost, as well as  resolution of complex solution structures}. In our notation, we denote the splitting SLDG method using $Q^k$ approximation spaces as $Q^k$ SLDG-split, and  denote the non-splitting SLDG method using $P^k$ approximation space as  $P^k$ SLDG. For the non-splitting SLDG method, when the quadratic-curves quadrilateral approximation is used to approximate upstream cells, such a method is denoted as $P^k$ SLDG-QC. For the splitting method, we always employ the second order Strang splitting as discussed in Section~\ref{sec: splitting}. The time step is set as
\begin{equation}
\Delta t = \frac{ CFL }{ \frac{a }{\Delta x} + \frac{b }{\Delta y} },
\end{equation}
where $a$ and $b$ are maximum transport speed in $x$- and $y$-directions, respectively. $CFL$ number is specified for comparison with the RKDG formulation, whose $CFL$ upper bound is $1/(2k+1)$ with $k$ being the  degree of the polynomial.

\subsection{2D linear passive-transport problems} \label{num_transport}

\begin{exa}
(Linear transport equation with constant coefficient.) Consider
\begin{equation}
u_t + u_x + u_y = 0, \ (x,y)\in[-\pi,\pi]^2
\label{linear2d}
\end{equation}
with  periodic boundary conditions and the initial condition
$
u(x,y,0) = \sin(x+y).
$
 The exact solution  is
$
u(x,y,t) = \sin(x+y-2t).
$
Note that for this example, since the convection operators in $x$- and $y$-directions commute, there is no splitting error.
Table \ref{table:2dlinear} summarizes the $L^2$ and $L^\infty$
errors, the associated orders of convergence, and CPU cost for $Q^k$  SLDG-split and $P^k$ SLDG for $k=1,2$.
We let the final time $T=\pi$, and consider  $CFL=2.5$ and $CFL=10.5$.
  The expected $k+1$-th order  convergence is observed for all cases.
It is also observed that the error magnitude by $Q^k$ SLDG-split  is smaller than that by $P^k$ SLDG. This is because $Q^k$ has more degrees of freedom than $P^k$ and hence enjoys better approximation property, even though their approximation rates are the same.
By comparing the CPU cost of $Q^k$  SLDG-split and $P^k$ SLDG with the same $CFL$, we observe that  $Q^k$  SLDG-split is more efficient than $P^k$ SLDG, which is ascribed to the absence of splitting error as well as the simplicity of the 1D implementation procedure for $Q^k$  SLDG-split.


\begin{table}[!ht] \footnotesize
\caption{$P^k$ SLDG and $Q^k$ SLDG-split $(k=1,2)$ for \eqref{linear2d} with $u(x,y,0)=\sin (x+y)$ at $T=\pi$.}
\vspace{0.1in}
\begin{tabular}{c cc cc cc cc cc}

\hline
 Mesh & \multicolumn{5}{c}{$CFL=2.5$}
     &\multicolumn{5}{c}{$CFL=10.5$ } \\
   \cmidrule(lr){2-6} \cmidrule(lr){7-11}
{   }  &{$L^2$ error} & Order  & {$L^{\infty}$ error} & Order & CPU &{$L^2$ error} & Order & {$L^{\infty}$ error} & Order & CPU \\
\hline
    \multicolumn{11}{l}{ $P^1$ SLDG}
     \\

    $20^2$ &     7.24E-03 &    --      &    4.07E-02 & -- & 0.02
     & 7.43E-03 &  --  &  4.33E-02   &   --  &  0.01 \\
    $40^2$ &     1.82E-03 &     1.99 &     1.04E-02 &     1.96 & 0.14
     &  1.82E-03 &     2.03 &     1.05E-02 &     2.04  & 0.03\\
    $80^2$ &     4.55E-04 &     2.00 &     2.63E-03 &     1.99 & 1.02
     &    4.55E-04 &     2.00 &     2.64E-03 &     1.99 &  0.31\\
    $160^2$ &     1.14E-04 &     2.00 &     6.60E-04 &     2.00 & 7.89
     &   1.14E-04 &     2.00 &     6.61E-04 &     2.00 & 2.14\\ \hline
\multicolumn{11}{l}{ $Q^1$ SLDG-split}
     \\

    $20^2$ &     1.91E-03 &    --      &   4.04E-03 & -- & 0.02
     &  2.21E-03 &  --  & 3.92E-03  &   --  &  0.01 \\
    $40^2$ &     4.71E-04 &     2.02 &     9.55E-04 &     2.08 & 0.13
     &   6.76E-04 &     1.71 &     1.35E-03 &     1.53   & 0.03\\
    $80^2$ &     1.17E-04 &     2.00 &     2.31E-04 &     2.05 & 0.84
     &    1.17E-04 &     2.53 &     2.24E-04 &     2.59  & 0.25\\
    $160^2$ &    2.93E-05 &     2.00 &     5.67E-05 &     2.03 & 6.31
     &    2.93E-05 &     2.00 &     5.59E-05 &     2.01 & 1.56\\ \hline
%
\hline
    \multicolumn{11}{l}{ $P^2$ SLDG}
     \\

    $20^2$ &     3.54E-04 &  --  &    2.34E-03 & -- & 0.05
     & 3.64E-04 &  -- &2.25E-03 & -- & 0.02\\
    $40^2$ &     4.42E-05 &     3.00 &     2.92E-04 &     3.01 & 0.30
     &4.41E-05 &     3.05 &     2.89E-04 &     2.96 &  0.09\\
    $80^2$ &    5.53E-06 &     3.00 &     3.63E-05 &     3.00 & 2.25
     &    5.52E-06 &     3.00 &     3.62E-05 &     3.00& 0.62\\
    $160^2$ &     6.91E-07 &     3.00 &     4.54E-06 &     3.00 & 17.59
     &    6.91E-07 &     3.00 &     4.54E-06 &     2.99  &  4.58\\
\hline
\multicolumn{11}{l}{ $Q^2$ SLDG-split }
     \\

    $20^2$ &     4.18E-05 &  --  &    1.08E-04 & -- & 0.05
     & 7.08E-05 &  -- &1.71E-04 & -- & 0.01\\
    $40^2$ &     5.43E-06 &     2.94 &     1.38E-05 &     2.96& 0.30
     &8.38E-06 &     3.08 &     2.28E-05 &     2.91 &  0.09\\
    $80^2$ &    6.79E-07 &     3.00 &     1.73E-06 &     3.00& 2.22
     &   6.60E-07 &     3.67 &     1.71E-06 &     3.73 & 0.61\\
    $160^2$ &    8.49E-08 &     3.00 &     2.16E-07 &     3.00& 17.83
     &    8.51E-08 &     2.96 &     2.17E-07 &     2.98 &  4.30\\
\hline
%
\end{tabular}
\label{table:2dlinear}
\end{table}

\end{exa}

\begin{exa}
(Rigid body rotation.) Consider
\begin{equation}
u_t - (yu)_x + (xu)_y =0, \  (x,y)\in[-2\pi,2\pi]^2.
\label{rigid1}
\end{equation}
We first consider a circular symmetry initial condition $u(x,y,0)=\exp(-x^2-y^2)$ and run the simulations to $T=20\pi$ (10 periods of rotation) with $CFL=2.5$ and $CFL=10.5$. Table \ref{table:rigid} summarizes the $L^2$ and $L^\infty$
errors, the associated orders of convergence and CPU cost of both SLDG methods.
It is observed that $Q^k$ SLDG-split produces results with smaller error magnitude and at the same time requires less CPU time compared with $P^k$ SLDG. Indeed, even with a large $CFL=10.5$, the spatial error from $Q^k$ SLDG-split discretization still dominates the splitting error due to the symmetry of the solution.
 Then, we take another initial condition $u(x,y,0)=\exp (-x^2-10y^2)$, for which the circular symmetry no longer holds. Table \ref{table:rigid1} summarizes the $L^2$ and $L^\infty$ errors, the associated orders of convergence and CPU cost of both methods  with $CFL=2.5$ and $CFL=10.5$ at $T = 20\pi$. We observe that, for $CFL=2.5$, $Q^k$ SLDG-split  performs better than  $P^k$ SLDG in terms of CPU efficiency, due to the fact that the the spatial error  still dominates the splitting error.  While for a large $CFL=10.5$, the splitting error in time becomes significant and hence dominates the total numerical error for $Q^2$ SLDG-split. In particular, second order  convergence due to Strang splitting is observed, and the error magnitude of $Q^2$ SLDG-split is much larger than that of  $P^k$ SLDG with the same configuration. To better understand the splitting error in time,  we test the temporal accuracy of both schemes by varying $CFL$  with a fixed  mesh of $160\times160$ cells, see Table \ref{table:rigid_time}. It is observed that $Q^2$ SLDG-split suffers a second order splitting error, and the approximation quality deteriorates quickly as $CFL$ increases. On the other hand, $P^2$ SLDG with $CFL$ as large as $25$ can still provide an accurate approximation to the solution. The error magnitude stays nearly the same as $CFL$ increases, indicating the dominating spatial error and relatively small temporal error even with large $CFLs$.

\begin{table}[!ht] \footnotesize
\caption{$P^k$ SLDG and $Q^k$ SLDG-split $(k=1,2)$ for \eqref{rigid1} with $u(x,y,0)=\exp (-x^2-y^2)$ at $T=20\pi$.}
\vspace{0.1in}
\begin{tabular}{c cc cc cc cc cc}

\hline
 Mesh & \multicolumn{5}{c}{$CFL=2.5$}
     &\multicolumn{5}{c}{$CFL=10.5$ } \\
   \cmidrule(lr){2-6} \cmidrule(lr){7-11}
{   }  &{$L^2$ error} & Order  & {$L^{\infty}$ error} & Order & CPU &{$L^2$ error} & Order & {$L^{\infty}$ error} & Order & CPU \\
\hline
    \multicolumn{11}{l}{ $P^1$ SLDG}
     \\

    $20^2$ &   4.60E-02 &    --      &  5.25E-01  & -- & 1.00
     & 2.88E-02 &  --  &  3.62E-01  &   --  & 0.34 \\
    $40^2$ &   1.50E-02 &     1.62 &     1.95E-01 &     1.43 & 7.73
     &  7.29E-03 &     1.98 &     1.20E-01 &     1.59 & 2.06\\
    $80^2$ &   2.70E-03 &     2.47 &     3.76E-02 &     2.38 & 59.30
     &  1.19E-03 &     2.62 &     2.71E-02 &     2.15 & 15.48 \\
    $160^2$ & 3.88E-04 &     2.80 &     6.85E-03 &     2.46  &  468.88
     &  1.82E-04 &     2.71 &     6.05E-03 &     2.17  & 120.94 \\ \hline
\multicolumn{11}{l}{ $Q^1$ SLDG-split}
     \\

    $20^2$ &     2.88E-02 &    --      &    2.12E-01 & -- & 0.88
     & 2.01E-02 &  --  & 1.78E-01   &   --  &  0.17 \\
    $40^2$ &     8.43E-03 &     1.77 &     9.61E-02 &     1.14 & 7.08
     &  5.72E-03 &     1.81 &     6.80E-02 &     1.39 &1.58\\
    $80^2$ &     1.27E-03 &     2.74 &     1.56E-02 &     2.62 & 55.94
     &     8.49E-04 &     2.75 &     1.30E-02 &     2.39 &  11.58\\
    $160^2$ &   1.77E-04 &     2.84 &     2.59E-03 &     2.60  & 452.00
     &    1.21E-04 &     2.81 &     2.43E-03 &     2.42 & 91.19\\

\hline\hline
    \multicolumn{11}{l}{ $P^2$ SLDG}
     \\

    $20^2$ &   5.86E-03   &  --  &   7.59E-02  & -- & 2.17
     & 2.13E-03 &  -- & 3.45E-02 & -- &  0.63 \\
    $40^2$ &   3.07E-04 &     4.25 &     5.84E-03 &     3.70 & 16.64
     & 1.53E-04 &     3.81 &     3.39E-03 &     3.35   &  4.45\\
    $80^2$ &  1.84E-05 &     4.06 &     4.23E-04 &     3.79 & 129.36
     &  1.50E-05 &     3.35 &     4.50E-04 &     2.91   & 33.55 \\
    $160^2$ & 1.92E-06 &     3.26 &     6.00E-05 &     2.82  & 1032.59
     & 1.72E-06 &     3.12 &     5.87E-05 &     2.94 & 262.78\\
\hline
\multicolumn{11}{l}{ $Q^2$ SLDG-split }
     \\

    $20^2$ &     2.01E-03  &  --  &    2.95E-02 & -- & 2.34
     & 2.13E-03 &  -- & 1.78E-02  & -- & 0.52\\
    $40^2$ &    1.32E-04 &     3.93 &     3.50E-03 &     3.08 & 17.94
     &1.29E-04 &     4.04 &     1.71E-03 &     3.38 &  4.03\\
    $80^2$ &   8.96E-06 &     3.88 &     2.83E-04 &     3.63 & 142.02
     &    7.54E-06 &     4.10 &     2.22E-04 &     2.95 & 31.92\\
    $160^2$ &   9.19E-07 &     3.29 &     1.90E-05 &     3.90  & 1116.86
     &    9.90E-07 &     2.93 &     1.73E-05 &     3.68 &    265.55 \\

  \hline
\end{tabular}
\label{table:rigid}
\end{table}

\begin{table}[!ht] \footnotesize
\caption{$P^k$ SLDG and $Q^k$ SLDG-split $(k=1,2)$ for \eqref{rigid1} with $u(x,y,0)=\exp (-x^2-10y^2)$ at $T=20\pi$.}
\vspace{0.1in}
\begin{tabular}{c cc cc cc cc cc}

\hline
 Mesh & \multicolumn{5}{c}{$CFL=2.5$}
     &\multicolumn{5}{c}{$CFL=10.5$ } \\
   \cmidrule(lr){2-6} \cmidrule(lr){7-11}
{   }  &{$L^2$ error} & Order  & {$L^{\infty}$ error} & Order & CPU &{$L^2$ error} & Order & {$L^{\infty}$ error} & Order & CPU \\
\hline
    \multicolumn{11}{l}{ $P^1$ SLDG}
     \\

    $20^2$ &   4.46E-02 &    --      &  8.11E-01  & -- & 1.05
     & 3.99E-02 &  --  &  7.25E-01  &   --  & 0.33 \\
    $40^2$ &   3.27E-02 &     0.45 &     5.77E-01 &     0.49 & 7.86
     &  2.69E-02 &     0.57 &     4.87E-01 &     0.57 & 2.05\\
    $80^2$ &   1.78E-02 &     0.88 &     2.87E-01 &     1.00 & 61.98
     &  1.27E-02 &     1.08 &     2.23E-01 &     1.13 & 15.42 \\
    $160^2$ & 5.93E-03 &     1.59 &     1.03E-01 &     1.48 &  478.59
     &  3.43E-03 &     1.89 &     6.65E-02 &     1.74  & 119.97 \\ \hline
\multicolumn{11}{l}{ $Q^1$ SLDG-split}
     \\

    $20^2$ &     3.14E-02 &    --      &     3.86E-01 & -- & 0.89
     & 3.46E-02 &  --  & 3.80E-01   &   --  &  0.14 \\
    $40^2$ &     2.62E-02 &     0.26 &     4.79E-01 &    -0.31 & 7.11
     &  2.42E-02 &     0.52 &     4.33E-01 &    -0.19 & 1.73\\
    $80^2$ &      1.39E-02 &     0.92 &     2.32E-01 &     1.05 & 53.70
     &    1.15E-02 &     1.07 &     2.03E-01 &     1.09  & 12.86\\
    $160^2$ &  4.17E-03 &     1.73 &     7.54E-02 &     1.62  & 422.23
     &   3.00E-03 &     1.94 &     5.77E-02 &     1.81 &  103.75\\

\hline\hline
    \multicolumn{11}{l}{ $P^2$ SLDG}
     \\

    $20^2$ &   2.73E-02   &  --  &   4.39E-01  & -- & 2.19
     & 3.49E-02 &  -- & 3.56E-01 & -- &  0.63 \\
    $40^2$ &   1.09E-02 &     1.33 &     1.75E-01 &     1.33  & 16.61
     & 7.06E-03 &     1.68 &     1.25E-01 &     1.50  &  4.45\\
    $80^2$ &  1.57E-03 &     2.80 &     2.79E-02 &     2.65  & 137.16
     & 7.09E-04 &     3.31 &     1.48E-02 &     3.08   & 33.59 \\
    $160^2$ & 7.93E-05 &     4.30 &     1.77E-03 &     3.98   & 1095.58
     & 3.17E-05 &     4.48 &     9.09E-04 &     4.02 & 261.95\\
\hline
\multicolumn{11}{l}{ $Q^2$ SLDG-split }
     \\

    $20^2$ &     2.23E-02  &  --  &    3.67E-01 & -- & 2.38
     & 3.63E-02  &  -- & 4.06E-01  & -- & 0.55\\
    $40^2$ &    7.10E-03 &     1.65 &     1.15E-01 &     1.67  & 18.47
     &1.44E-02 &     1.33 &     2.01E-01 &     1.02 &  4.44\\
    $80^2$ &    7.69E-04 &     3.21 &     1.37E-02 &     3.07 & 146.53
     &  3.62E-03 &     1.99 &     5.00E-02 &     2.01 & 35.33\\
    $160^2$ &    6.09E-05 &     3.66 &     1.13E-03 &     3.60  & 1107.16
     &   9.00E-04 &     2.01 &     1.20E-02 &     2.06 &  274.30\\

  \hline
\end{tabular}
\label{table:rigid1}
\end{table}

\begin{table}[!ht] 
\small
\caption{$P^2$ SLDG and $Q^2$ SLDG-split  for \eqref{rigid1} with $u(x,y,0)=\exp (-x^2-10y^2)$ at $T=20\pi$. A mesh of $160\times160$ cells is used.}
\vspace{0.1in}
\centering
\begin{tabular}{c cc cc c }

\hline
{  $CFL$ }  &{$L^2$ error} & Order  & {$L^{\infty}$ error} & Order & CPU  \\
\hline
%
    \multicolumn{6}{l}{ $P^2$ SLDG}
     \\

    5 &   5.62E-05 & -- &     1.33E-03  & -- & 487.08
       \\
    10 &   3.30E-05 &    -0.77 &     9.15E-04 &    -0.55 & 251.50
      \\
    15 &  2.45E-05 &    -0.73 &     8.63E-04 &    -0.14 &165.55
       \\
    20 &  2.09E-05 &    -0.55 &     8.16E-04 &    -0.20  & 129.88
      \\
      25 &  1.95E-05 &    -0.31 &     8.56E-04 &     0.22  & 109.30
      \\
\hline
\multicolumn{6}{l}{ $Q^2$ SLDG-split }
     \\

    5 &     2.06E-04  &  --  &    3.01E-03 & -- & 542.84
      \\
    10 &    8.16E-04 &     1.99 &     1.09E-02 &     1.86 & 271.59
     \\
    15 &    1.84E-03 &     2.00 &     2.43E-02 &     1.98 & 184.47
     \\
    20 &    3.26E-03 &     2.00 &     4.27E-02 &     1.96  & 137.52
      \\
    25 &    5.11E-03 &     2.01 &     6.70E-02 &     2.02  & 109.41
     \\
\hline
%
\end{tabular}
\label{table:rigid_time}
\end{table}

\end{exa}

\begin{exa}
(Swirling deformation flow.)
We consider solving
\begin{equation}
u_t - \left( \cos^2 \left(\frac{x}{2} \right)\sin(y) g(t) u \right)_x
+ \left(  \sin(x) \cos^2 \left(\frac{y}{2}\right) g(t) u \right)_y =0,
\
(x,y)\in[-\pi,\pi]^2,
\label{swirling}
\end{equation}
 where $g(t) = \cos\left( \frac{\pi t}{T} \right)\pi$.
 $T=1.5$.
 The initial condition is set
 to be the following smooth cosine bells (with $C^5$ smoothness),
\begin{equation}
u(x,y,0) =
\begin{cases}
r_0^b \cos^6 \left(  \frac{r^b}{2r_0^b} \pi \right), & \text{if}\  r^b <r_0^b,\\
0,  & \text{ otherwise},
\end{cases}
\label{cosine_bell}
\end{equation}
where $r_0^b = 0.3\pi$, and $r^b=\sqrt{  (x-x_0^b)^2 + (y-y_0^b)^2 }$ denotes the distance between $(x,y)$ and the center of the cosine bell $(x_0^b, y_0^b) = ( 0.3\pi ,0 )$.
%
Note that this swirling deformation flow  introduced in \cite{leveque1996high}  is more challenging than the rigid body rotation due to the space- and time-dependent flow field.
In particular, along the direction of the flow, the initial function becomes largely deformed at $t=T/2$, then goes back to its initial shape at $t = T$ as the flow reverses. If this problem is solved up to $T$, we call such a procedure  one full evolution. Note that unlike the previous two examples, the upstream cells are no longer rectangular, and hence the quadratic-curved quadrilateral approximation is  expected to capture the geometry of the upstream  cells more accurately than the quadrilateral approximation. However, regarding the efficiency,  the quadrilateral approximation is already adequate for $k=1$, since the DG appropriation is  second accurate. The efficiency benefit of using more costly quadratic-curved quadrilateral approximation is more evident for $k=2$.

We test order of accuracy for $P^k$ SLDG(-QC) and $Q^k$ SLDG-split for $k=1,2$ up to $T/2$, i.e.,  half of one full evolution, and summarize the results in Tables \ref{table:swirling_ref}. Note that as the exact solution is not available, we pre-compute a reference solution by $P^2$ SLDG-QC with a refined mesh of $320\times320$ cells and $CFL=2.5$.
As expected, the $(k+1)^{th}$ order convergence is observed for $P^k$ SLDG(-QC) and the second order convergence is observed due to the splitting error. More specifically, when comparing $P^1$ SLDG and  $Q^1$ SLDG-split, error magnitude from the splitting one is relatively smaller, which again is ascribed  to the fact that despite the same convergence rate, $Q^1$  delivers better approximation performance than $P^1$. Further, when comparing $P^2$ SLDG-QC and  $Q^2$ SLDG-split, the splitting method does not perform as well as the non-splitting one in terms of error magnitude, due to its splitting error.

We then perform the accuracy test in time with a fixed mesh of $160\times160$ cells and different $CFL$s ranging from $5$ to $25$. The reference solution is computed by $P^2$ SLDG-QC with the same mesh but using a small $CFL=0.1$. We report the results in Table \ref{swirling_time1}. The second order Strang splitting error is clearly observed in the Table. It is also observed that,  $P^2$ SLDG-QC significantly outperforms $Q^2$ SLDG-split in terms of error magnitude. The two methods consume comparable CPU time with the same simulation configuration.

We also test the convergence for {\em one full} evolution and summarize the results in Table \ref{table:swirling} for $P^k$ SLDG(-QC) and $Q^k$ SLDG-split, $k=1,2$, with $CFL=2.5$ and $CFL=10.5$. It is observed that the errors after one full evolution are less than those at one half  evolution. Such a super performance is ascribed to the special error cancellations (coming from some symmetry in the solution movement along the flow) as explained in \cite{blossey2008selective}. In addition, the error cancellation also makes the splitting error less pronounced for $Q^k$ SLDG-split. This is an example, in which special cancellation of errors plays a role in comparing schemes' performances.

\begin{table}[!ht] \footnotesize
\caption{Swirling deformation flow. $P^k$ SLDG(-QC) and $Q^k$ SLDG-split ($k=1,2$) for \eqref{swirling} with the smooth cosine bells \eqref{cosine_bell} at $T/2=0.75$. A reference solution is solved by $P^2$ SLDG-QC with $CFL=2.5$ and  a mesh of $320\times320$ cells. }
\vspace{0.1in}
\begin{tabular}{c cc cc cc cc cc}
\hline
 Mesh & \multicolumn{5}{c}{$CFL=2.5$}
     &\multicolumn{5}{c}{$CFL=10.5$ } \\
   \cmidrule(lr){2-6} \cmidrule(lr){7-11}
{   }  &{$L^2$ error} & Order  & {$L^{\infty}$ error} & Order & CPU &{$L^2$ error} & Order & {$L^{\infty}$ error} & Order & CPU \\
\hline
    \multicolumn{11}{l}{ $P^1$ SLDG}
     \\
    $20^2$ &  1.08E-02  &    --      & 2.49E-01  & -- & 0.03
     &  9.27E-03  &  --  &   2.34E-01 &   --  & 0.02 \\
    $40^2$ &  2.98E-03 &     1.86 &     1.04E-01 &     1.26  & 0.25
     &   3.50E-03 &     1.40 &     1.31E-01 &     0.83 & 0.09 \\
    $80^2$ &  7.88E-04 &     1.92 &     3.23E-02 &     1.69  & 1.70
     &  8.47E-04 &     2.05 &     3.78E-02 &     1.80  & 0.53\\
    $160^2$ &  2.14E-04 &     1.88 &     9.02E-03 &     1.84 &  12.98
     &  2.30E-04 &     1.88 &     1.28E-02 &     1.56  & 3.81 \\ \hline
\multicolumn{11}{l}{ $Q^1$ SLDG-split}
     \\
    $20^2$ &  6.68E-03  &    --      &  1.42E-01  & -- & 0.03
     &  3.18E-02   &  --  &     4.19E-01   &   --  & 0.02 \\
    $40^2$ &  1.66E-03 &     2.01 &     4.37E-02 &     1.70 & 0.23
     &  6.27E-03 &     2.34 &     9.15E-02 &     2.20 & 0.06 \\
    $80^2$ &  4.32E-04 &     1.95 &     1.10E-02 &     1.99 & 1.75
     &   1.48E-03 &     2.08 &     2.27E-02 &     2.01 & 0.44 \\
    $160^2$ &  1.13E-04 &     1.94 &     2.68E-03 &     2.04  & 13.88
     &  3.69E-04 &     2.01 &     6.00E-03 &     1.92  & 3.48 \\ \hline
\hline
    \multicolumn{11}{l}{ $P^2$ SLDG-QC}
     \\

    $20^2$ &   2.43E-03  &  --  &   9.30E-02 & -- & 0.06
     & 2.76E-03  &  -- & 9.55E-02 & -- &  0.02 \\
    $40^2$ &   3.44E-04 &     2.82 &     2.13E-02 &     2.13  & 0.42
     & 4.54E-04 &     2.60 &     1.71E-02 &     2.48 &  0.14\\
    $80^2$ &  4.56E-05 &     2.91 &     3.08E-03 &     2.79  &3.27
     &  5.52E-05 &     3.04 &     2.41E-03 &     2.83  & 0.86\\
    $160^2$ &  5.80E-06 &     2.98 &     3.70E-04 &     3.06   &25.58
     & 6.21E-06 &     3.15 &     3.04E-04 &     2.99 & 6.36\\
\hline
\multicolumn{11}{l}{ $Q^2$ SLDG-split }
     \\
    $20^2$ &     1.56E-03  &  --  &    2.60E-02 & -- & 0.06
     &  3.32E-02  &  -- &  4.36E-01  & -- & 0.02\\
    $40^2$ &    3.53E-04 &     2.14 &     6.47E-03 &     2.01 & 0.45
     & 6.26E-03 &     2.41 &     8.89E-02 &     2.30  &  0.11\\
    $80^2$ &    8.36E-05 &     2.08 &     1.43E-03 &     2.18 & 3.69
     &   1.47E-03 &     2.09 &     2.12E-02 &     2.06  & 0.92\\
    $160^2$ &   2.05E-05 &     2.03 &     3.29E-04 &     2.11  & 29.17
     &  3.60E-04 &     2.03 &     5.18E-03 &     2.03  &  7.30\\
  \hline
\end{tabular}
\label{table:swirling_ref}
\end{table}

\begin{table}[!ht] 
\small
\caption{Swirling deformation flow. Order of accuracy in time for $P^2$ SLDG-QC and $Q^2$ SLDG-split  for \eqref{swirling} with the smooth cosine bells \eqref{cosine_bell} at $T/2=0.75$  by comparing numerical solutions with
a reference solution from the corresponding scheme with $CFL=0.1$.   A mesh of $160\times160$ is used. }
\vspace{0.1in}
\centering
\begin{tabular}{c cc cc c }

\hline
{  $CFL$ }  &{$L^2$ error} & Order  & {$L^{\infty}$ error} & Order & CPU  \\
\hline

    \multicolumn{6}{l}{ $P^2$ SLDG-QC}
     \\

    5 &   1.85E-06 & -- &     7.34E-05 &      --         &8.48
       \\
    10 &   4.13E-06 &     1.16 &     2.16E-04 &     1.56  & 4.38
      \\
    15 &   7.22E-06 &     1.38 &     4.74E-04 &     1.94 & 3.05
       \\
    20 &  1.13E-05 &     1.55 &     8.10E-04 &     1.86 &2.34
      \\
    25 &   1.36E-05 &     0.84 &     8.33E-04 &     0.13 &1.98
      \\
\hline
\multicolumn{6}{l}{ $Q^2$ SLDG-split }
     \\
    5 &   8.09E-05 & -- &     1.18E-03 &      --         &9.45
       \\
    10 &   3.23E-04 &     2.00 &     4.67E-03 &     1.99 &4.75
      \\
    15 &   7.28E-04 &     2.00 &     1.06E-02 &     2.01 & 3.17
       \\
    20 &  1.30E-03 &     2.01 &     1.88E-02 &     2.01 & 2.36
      \\
    25 &  2.08E-03 &     2.11 &     3.02E-02 &     2.11 & 2.00
      \\
\hline

\end{tabular}
\label{swirling_time1}
\end{table}

\begin{table}[!ht] \footnotesize
\caption{Swirling deformation flow. $P^k$ SLDG(-QC) and $Q^k$ SLDG-split ($k=1,2$) for \eqref{swirling} with the smooth cosine bells \eqref{cosine_bell} at $T=1.5$.}
\vspace{0.1in}
\begin{tabular}{c cc cc cc cc cc}
\hline
 Mesh & \multicolumn{5}{c}{$CFL=2.5$}
     &\multicolumn{5}{c}{$CFL=10.5$ } \\
   \cmidrule(lr){2-6} \cmidrule(lr){7-11}
{   }  &{$L^2$ error} & Order  & {$L^{\infty}$ error} & Order & CPU &{$L^2$ error} & Order & {$L^{\infty}$ error} & Order & CPU \\
\hline
    \multicolumn{11}{l}{ $P^1$ SLDG}
     \\
    $20^2$ &  1.41E-02  &    --      & 2.94E-01  & -- & 0.08
     &  1.06E-02  &  --  &   2.11E-01 &   --  & 0.01 \\
    $40^2$ &  3.85E-03 &     1.87 &     9.33E-02 &     1.65  & 0.44
     &  5.47E-03 &     0.95 &     2.03E-01 &     0.05 & 0.13\\
    $80^2$ &  9.55E-04 &     2.01 &     3.12E-02 &     1.58 & 3.55
     &  1.57E-03 &     1.80 &     9.23E-02 &     1.14 & 0.89\\
    $160^2$ &  2.60E-04 &     1.87 &     1.33E-02 &     1.23 &  27.58
     &  6.13E-04 &     1.36 &     3.15E-02 &     1.55  & 6.61 \\ \hline
\multicolumn{11}{l}{ $Q^1$ SLDG-split}
     \\
    $20^2$ &  9.88E-03  &    --      &  1.87E-01  & -- & 0.06
     &  1.10E-02   &  --  &    1.55E-01  &   --  & 0.02 \\
    $40^2$ &  2.15E-03 &     2.20 &     5.41E-02 &     1.79 & 0.50
     &  2.73E-03 &     2.02 &     4.75E-02 &     1.71 & 0.14\\
    $80^2$ &  4.03E-04 &     2.41 &     1.28E-02 &     2.08 &4.04
     &   3.86E-04 &     2.82 &     7.45E-03 &     2.67 & 1.00\\
    $160^2$ &  8.47E-05 &     2.25 &     2.75E-03 &     2.22  &  31.63
     &  5.50E-05 &     2.81 &     1.41E-03 &     2.40  & 7.56 \\ \hline
\hline
    \multicolumn{11}{l}{ $P^2$ SLDG-QC}
     \\

    $20^2$ &   2.68E-03  &  --  &   5.78E-02 & -- & 0.14
     & 3.38E-03  &  -- & 1.11E-01 & -- &  3.13 \\
    $40^2$ &   3.72E-04 &     2.85 &     1.23E-02 &     2.23  & 0.83
     & 5.50E-04 &     2.62 &     1.64E-02 &     2.75 &  0.22\\
    $80^2$ &  5.69E-05 &     2.71 &     3.12E-03 &     1.98  &6.38
     &  7.99E-05 &     2.78 &     4.50E-03 &     1.87  & 1.64 \\
    $160^2$ &  8.56E-06 &     2.73 &     6.63E-04 &     2.23   &50.78
     & 1.21E-05 &     2.73 &     6.46E-04 &     2.80 & 12.30\\
\hline
\multicolumn{11}{l}{ $Q^2$ SLDG-split }
     \\
    $20^2$ &     1.13E-03  &  --  &    1.97E-02 & -- & 0.14
     & 1.16E-02  &  -- & 1.74E-01  & -- & 0.03\\
    $40^2$ &    1.26E-04 &     3.17 &     2.85E-03 &     2.79 & 1.05
     &2.55E-03 &     2.19 &     3.57E-02 &     2.28 &  0.25\\
    $80^2$ &    1.49E-05 &     3.08 &     4.15E-04 &     2.78 & 8.19
     &   2.95E-04 &     3.11 &     4.19E-03 &     3.09  & 2.03\\
    $160^2$ &    1.86E-06 &     3.01 &     4.91E-05 &     3.08  & 65.86
     &  2.39E-05 &     3.63 &     3.84E-04 &     3.45  &  15.53\\
  \hline
\end{tabular}
\label{table:swirling}
\end{table}

\end{exa}

\subsection{Vlasov-Poisson system} \label{num_VP}

In this subsection, we  compare the performances of the splitting and non-splitting SLDG methods for solving the 1D1V VP system \eqref{equation:vp}. 	We focus on the strong Landau damping test. The observations for other benchmark tests including the two stream instability I \cite{FilbetS}, two stream instability II \cite{umeda2008conservative} and bump-on-tail instability \cite{arber2002critical,gucclu2014arbitrarily} are similar, thus are omitted for brevity. For the non-splitting SLDG method, we employ high order characteristics tracing schemes as proposed in \cite{cai2018high}. In particular, we couple $P^1$ SLDG with a second order characteristic tracing scheme and couple $P^2$ SLDG-QC with a third order characteristic tracing scheme, and  denote the resulting methods as $P^1$ SLDG-time2 and $P^2$ SLDG-QC-time3, respectively. We apply the positivity-preserving limiter \cite{zhangshu2010} for all nonlinear test examples.

%
%
%

\begin{exa}
(Strong Landau damping.) Consider the strong Landau damping for the VP system with
the initial condition being a perturbed Maxwellian equilibrium
\begin{equation}
f(x,v,t=0) = \frac{ 1 }{\sqrt{2\pi} } ( 1+ \alpha \cos(k_0x)  ) \exp \left( -\frac{v^2}{2} \right) ,
\label{init}
\end{equation}
where $\alpha = 0.5$ and $k_0=0.5$.
This problem has been numerically investigated by several authors, e.g. see \cite{qiu2011positivity,crouseilles2011discontinuous,xiong2014high,zhu2016h,huang2016semi}.
First, we perform the comparison for $Q^k$ SLDG-split and $P^k$ SLDG-(QC)-time$k+1$, $k=1,2$. The results are reported in Table \ref{temporal_1}.
In the simulation, we use a fixed mesh of $200\times 200$ cells. We let $T=2$ and  $CFL$ range from $5$ to $25$.
The errors are computed by comparing numerical solutions with the reference solution by $P^2$ SLDG-QC-time3 with the same mesh but using a small $CFL=0.01$. Expected orders of convergence are observed. By taking a closer look at the error table, we observe that the error magnitude from both methods is comparable, while the splitting method takes less CPU time. The methods also perform similarly in conserving physical invariants and in resolving solution structures, see Figure \ref{landau_norms} and Figure \ref{landau_40}, respectively. The splitting scheme outperforms the non-splitting one, as the temporal splitting error is not dominating, and $Q^k$ provides better approximations than $P^k$ polynomials.

\begin{table}[!ht]\small
\caption{Strong Landau damping.
 Temporal order of convergence and CPU comparison between $Q^k$ SLDG-split and $P^k$ SLDG-(QC)-time$k+1$ by comparing
numerical solutions with a reference solution from the corresponding scheme with $CFL=0.01$.  $T=2$. The mesh of $200\times200$ cells is used.
 }
\vspace{0.1in}
\centering
\begin{tabular}{c cc  cc c  cc cc c}
\hline

{$CFL$}   &{$L^2$ error} & Order    &{$L^\infty$ error} & Order &CPU
 &{$L^2$ error} & Order    &{$L^\infty$ error} & Order &CPU\\
\hline
& \multicolumn{5}{c}{$Q^1$ SLDG-split}   & \multicolumn{5}{c}{$P^1$ SLDG-time2} \\
 \cmidrule(lr){2-6} \cmidrule(lr){7-11}
  5  &   1.70E-05 & &     1.10E-04 &         & 4.70 &   1.44E-04 & &     7.09E-04 &  &  10.49\\
  10  &  4.67E-05 &     1.46 &     2.71E-04 &     1.31 & 2.39 &   5.63E-04 &     1.97 &     2.56E-03 &     1.85  & 6.29 \\
  15  &  9.03E-05 &     1.63 &     4.01E-04 &     0.96 & 1.64 &   1.26E-03 &     1.98 &     5.60E-03 &     1.93  & 3.74\\
  20 &  1.61E-04 &     2.01 &     7.53E-04 &     2.19 & 1.21 &   2.24E-03 &     2.01 &     1.01E-02 &     2.04 & 2.77\\
  25 &  2.51E-04 &     2.00 &     1.04E-03 &     1.47 & 1.02 &   3.44E-03 &     1.93 &     1.54E-02 &     1.91 & 2.49\\
\hline

& \multicolumn{5}{c}{$Q^2$ SLDG-split}   & \multicolumn{5}{c}{$P^2$ SLDG-QC-time3} \\
 \cmidrule(lr){2-6} \cmidrule(lr){7-11}
  5  &  9.40E-06 & &     3.67E-05 &         & 16.29 &   5.07E-06 & &     5.40E-05 &  & 28.10\\
  10  &  3.80E-05 &     2.01 &     1.46E-04 &     1.99  & 8.16 &   1.83E-05 &     1.85 &     1.66E-04 &   1.62  & 11.40 \\
  15  &  8.70E-05 &     2.04 &     3.31E-04 &     2.03  & 5.60 &  5.54E-05 &     2.73 &     3.39E-04 &     1.76 & 7.85 \\
  20 &  1.44E-04 &     1.76 &     5.54E-04 &     1.78 & 4.11 &   1.32E-04 &     3.03 &     7.39E-04 &     2.71 & 5.89 \\
  25 &  2.41E-04 &     2.30 &     9.10E-04 &     2.22 & 3.40 &   2.55E-04 &     2.94 &     1.38E-03 &     2.78 & 4.77\\
\hline
\end{tabular}
\label{temporal_1}
\end{table}


\begin{figure}[h!]
\centering                              
\includegraphics[height=70mm]{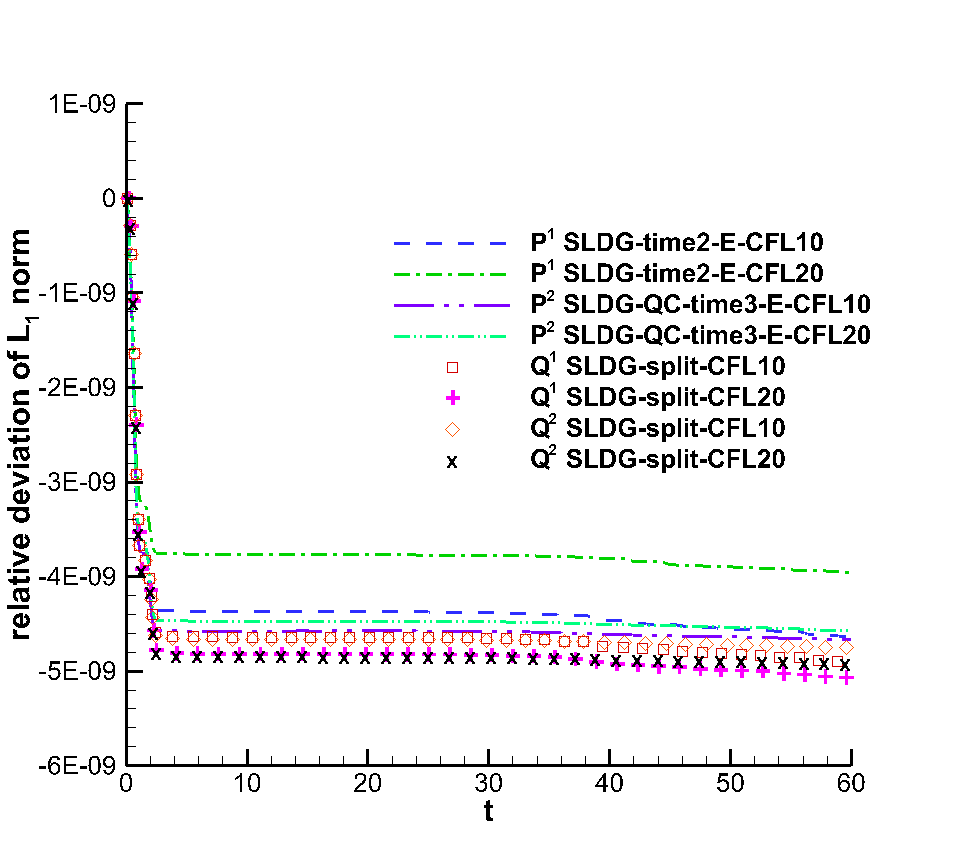}
\includegraphics[height=70mm]{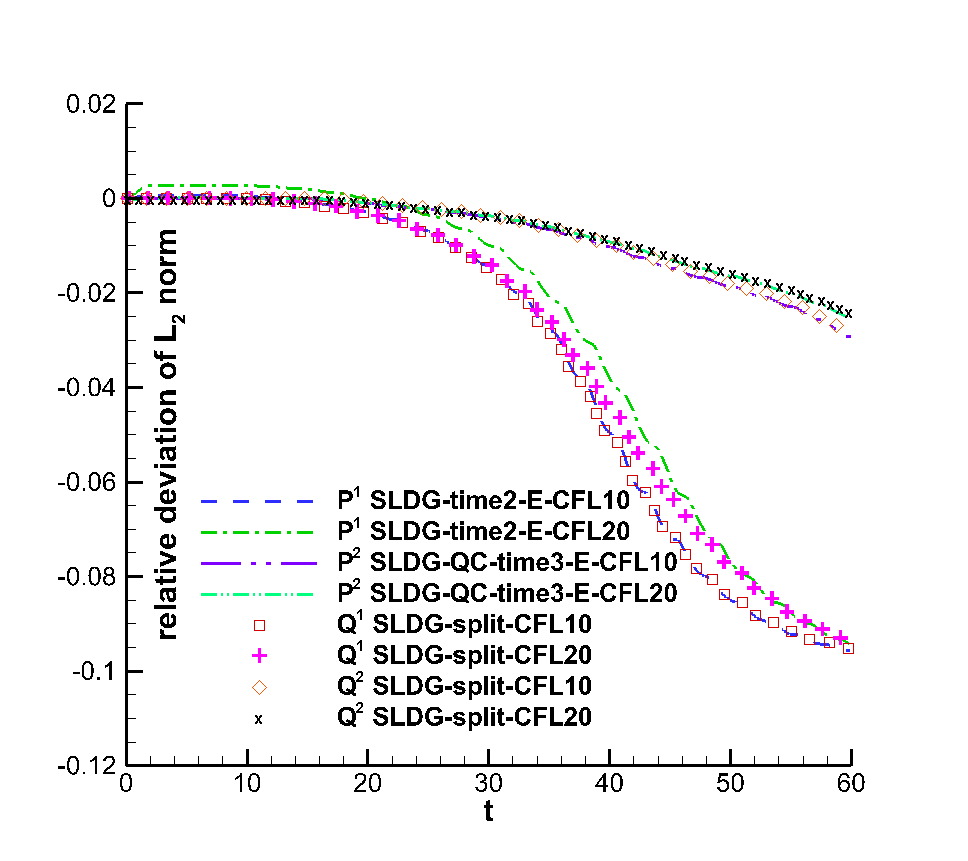}
\includegraphics[height=70mm]{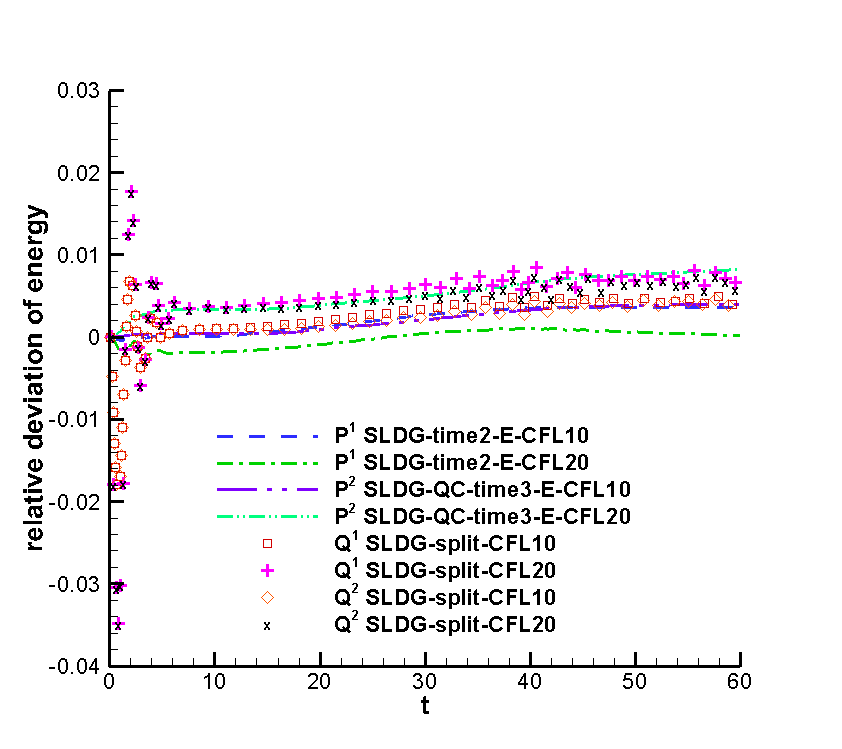}
\includegraphics[height=70mm]{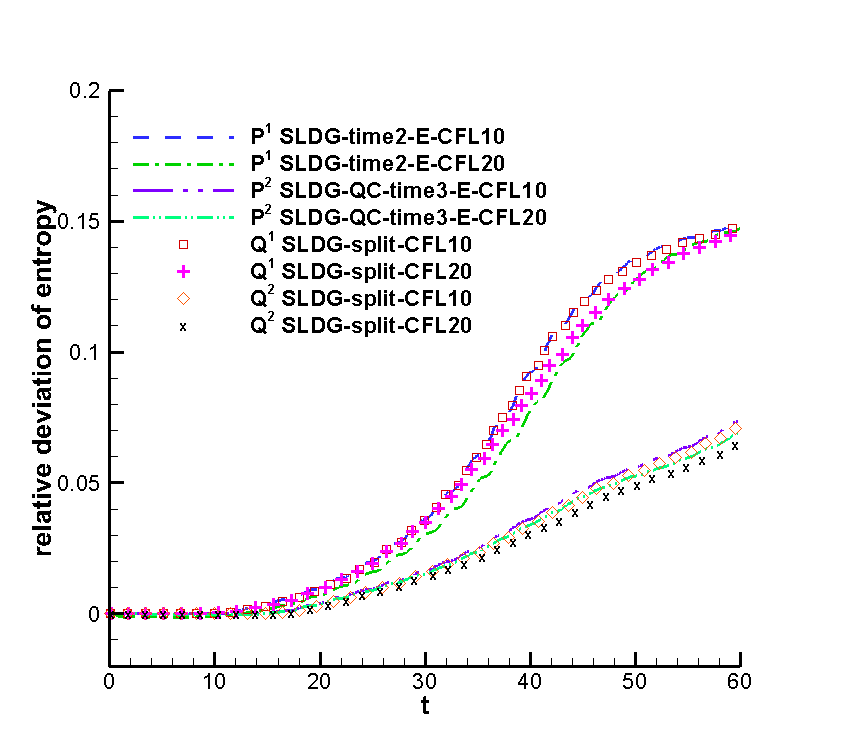}
\caption{Strong Landau damping. Time evolution of the relative deviations of $L^1$ (upper left) and $L^2$ (upper right) norms of the solution as well as the discrete kinetic energy (lower left) and entropy (lower right).}
\label{landau_norms}
\end{figure}

\begin{figure}[h!]
\centering                              
\includegraphics[height=70mm]{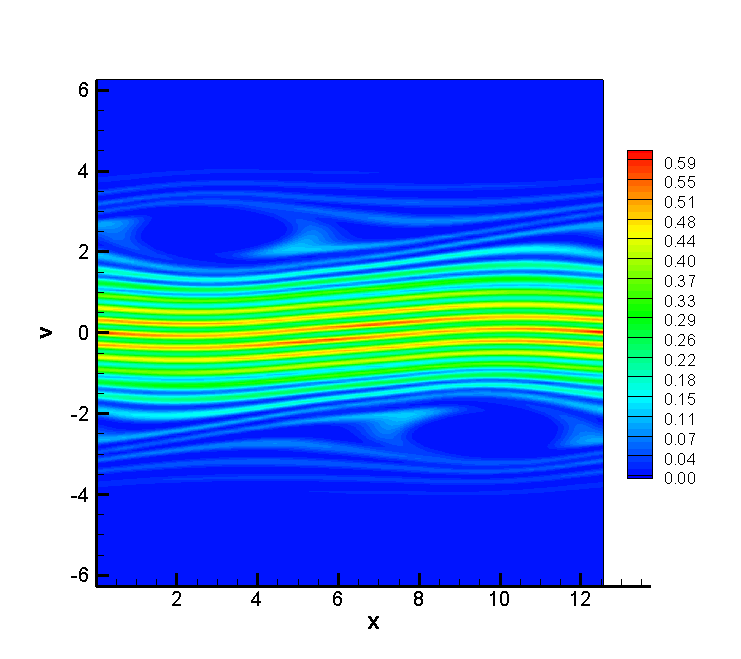}
\includegraphics[height=70mm]{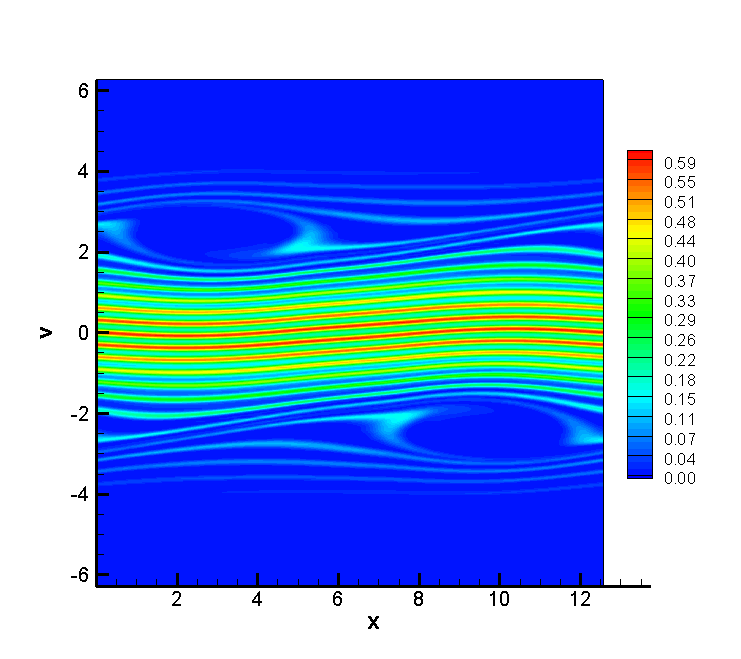}
\includegraphics[height=70mm]{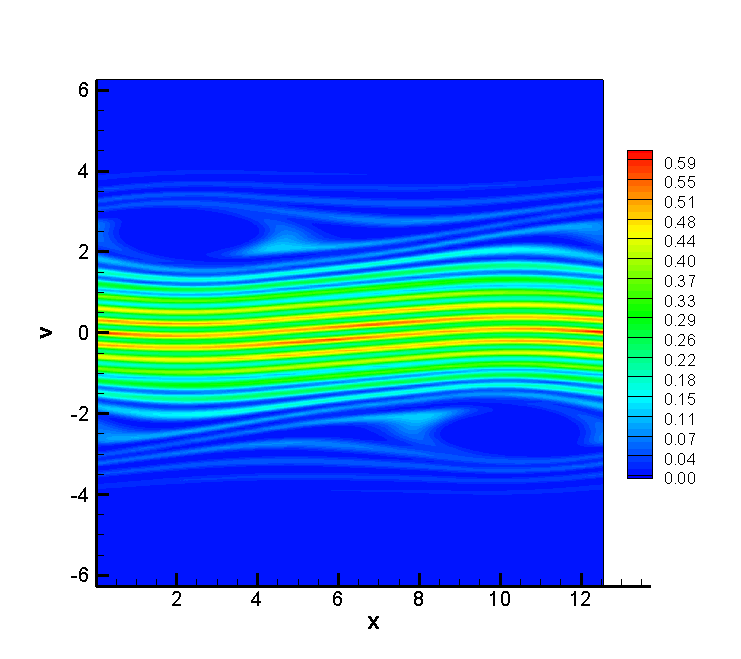}
\includegraphics[height=70mm]{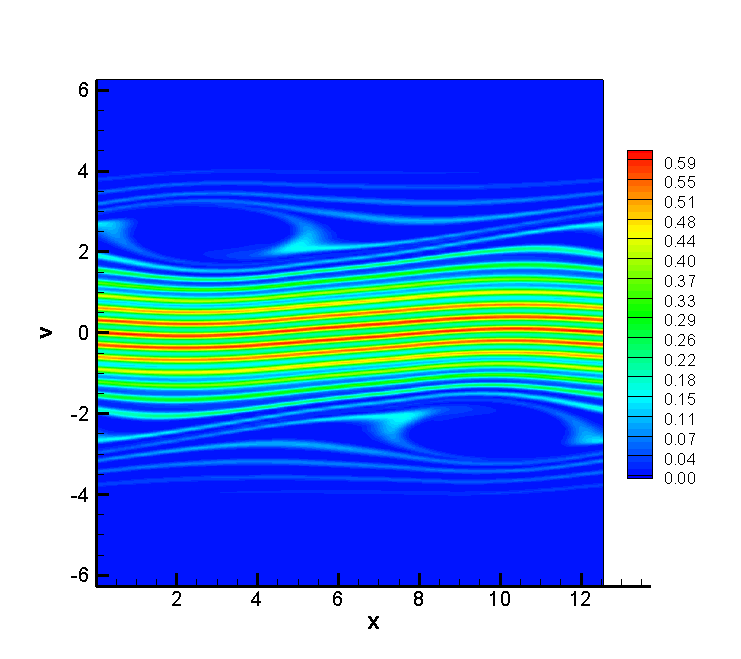}
\caption{Strang Landau damping with the spatial mesh of $160\times160$. $T=40$. $CFL=10$.
Upper left: $Q^1$ SLDG-split. Upper right: $Q^2$ SLDG-split. Bottom left: $P^1$ SLDG-time2.
Bottom right: $P^2$ SLDG-QC-time3.}
\label{landau_40}
\end{figure}

\end{exa}

\subsection{2D incompressible Euler equations and the guiding center Vlasov model} \label{num_Euler}
In this subsection, we compare the performance of the SLDG methods for solving the 2D incompressible Euler equation in the vorticity-stream function formulation \eqref{Euler} and the guiding center Vlasov model \eqref{guiding_center}.
As with the VP system, we have to couple a high order characteristics tracing mechanism for the non-splitting SLDG method. For the notational convenience, below $P^k$ SLDG(-QC)-$P^{r}$ LDG-time$(k+1)$-CFL$s$ denotes the non-splitting SLDG scheme with $P^k$ polynomial space, using the $P^{r}$ LDG scheme for solving Poisson's equation,  the $(k+1)^{th}$ order scheme in characteristics tracing, and $CFL=s$. As mentioned in \cite{cai2018Euler}, the choice of using the $P^{k}$ or $P^{k+1}$ LDG solver, i.e., $r=k$ or $k+1$, for Poisson's equation is a trade-off between accuracy (effectiveness in resolving solutions) and CPU cost. In the simulations, we let $r=k+1$ for the accuracy test only and let $r=k$ otherwise. See \cite{cai2018Euler} for more detailed discussion.
\begin{exa}
(Accuracy and convergence test.) Consider the incompressible Euler equation \eqref{Euler} on the domain $[0,2\pi]\times[0,2\pi]$ with the initial condition
\begin{equation}
\omega(x,y,0) = -2\sin(x) \sin(y)
\end{equation}
and periodic boundary conditions. The exact solution stays stationary as $\omega(x,y,t) = -2\sin(x)\sin(y)$. Here we solve the problem up to $T=1$ with $CFL=1$. We test the accuracy and CPU cost for $Q^k$ SLDG-split-$P^{k+1}$ LDG-CFL1 and $P^k$ SLDG(-QC)-$P^{k+1}$ LDG-time$(k+1)$-CFL1, $k=1,2$, and summarize results in Table \ref{Euler_order}. Expected second and third order convergence is observed for $P^1$ SLDG-$P^2$ LDG-time2 and $P^2$ SLDG-QC-$P^3$ LDG-time3, respectively. Only first order of convergence is observed for the splitting scheme due to the first order error in tracing characteristics, even if the second order Strang splitting is used. It is evident that the non-splitting method is more efficient due to its genuine high order accuracy in both space and time.

 \begin{table}[!ht]\small
\caption{The incompressible Euler equations. Errors, orders and CPU times (sec) comparison between $Q^k$ SLDG-split and $P^k$ SLDG(-QC) $k=1,2$. $CFL=1$. $T=1$.}
\vspace{0.1in}
\centering
\begin{tabular}{c cc  cc c  cc cc c}
\hline

{ Mesh}   &{$L^2$ error} & Order    &{$L^\infty$ error} & Order &CPU
 &{$L^2$ error} & Order    &{$L^\infty$ error} & Order &CPU\\
\hline
& \multicolumn{5}{c}{$Q^1$ SLDG-split+$P^2$ LDG}   & \multicolumn{5}{c}{$P^1$ SLDG+$P^2$ LDG+time2} \\
 \cmidrule(lr){2-6} \cmidrule(lr){7-11}
  $20^2$  &  4.06E-02 & &     1.30E-01 &         & 0.23 &  1.57E-02 & &     8.55E-02 &  & 0.21\\
  $40^2$  &  2.00E-02 &     1.03 &     6.18E-02 &     1.07 & 2.97 & 8.55E-02&   1.97& 2.48E-02 &  1.78 & 2.47 \\
  $60^2$  &  1.34E-02 &     0.98 &     4.09E-02 &     1.02  & 16.55 &2.48E-02 &  1.98 &   1.16E-02 &  1.87& 12.68 \\
  $80^2$ &  1.00E-02 &     1.02 &     3.03E-02 &     1.05 & 75.91 & 1.02E-03 &   1.96 & 6.71E-03 &  1.91 & 49.12 \\
  $100^2$  &  8.03E-03 &     0.98 &     2.42E-02 &     1.00& 193.45 & 6.49E-04&  2.01 & 4.34E-03 &  1.95 & 131.11\\
\hline

& \multicolumn{5}{c}{$Q^2$ SLDG-split+$P^3$ LDG}   & \multicolumn{5}{c}{$P^2$ SLDG-QC+$P^3$ LDG+time3} \\
 \cmidrule(lr){2-6} \cmidrule(lr){7-11}
  $20^2$  &  4.04E-02 & &     1.17E-01 &         & 0.67 &  2.82E-03 & &   1.36E-02 &  & 1.63\\
  $40^2$  &  1.99E-02 &   1.02 &     5.80E-02 &     1.02  & 8.82& 3.56E-04 &  2.99 &   1.93E-03 &  2.81 & 17.63 \\
  $60^2$  &  1.34E-02 &   0.98 &     3.91E-02 &     0.97  &59.58 &1.04E-04 &  3.02 &  5.33E-04 &   3.18 & 113.77\\
  $80^2$ &  1.00E-02 &     1.02 &     2.91E-02 &     1.02 &205.18& 4.44E-05&  2.97 &   2.28E-04 &  2.96 & 338.95 \\
  $100^2$  & 8.03E-03 &     0.98 &     2.34E-02 &     0.98&667.46 & 2.24E-05 &  3.08 & 1.14E-04 &   3.09& 1003.01\\
\hline
\end{tabular}
\label{Euler_order}
\end{table}

\end{exa}

\begin{exa}
(The double shear layer problem \cite{bell1989second,minion1997performance,zhangshu2010}.)
We consider the model problem \eqref{Euler} in the domain $[0,2\pi]\times[0,2\pi]$, with periodic boundary conditions and
the initial condition
\begin{equation}
\omega(x,y,0)
=
\begin{cases}
\delta \cos(x) - \frac{1}{\rho} sech^2\left( \frac{y-\pi/2}{\rho} \right), & \text{if} \ y\leq \pi, \\
\delta \cos(x) + \frac{1}{\rho} sech^2\left( \frac{3\pi/2-y}{\rho} \right), & \text{if} \ y>\pi,
\end{cases}
\end{equation}
where $\delta =0.05$ and $\rho = \pi/15$.
As time evolves,
the shear layers quickly develop into roll-ups with thinner and thinner scales.
Eventually, with a fixed mesh, the full resolution of the structures will be lost for any methods. This problem has been tested by the high order Eulerian finite difference ENO/WENO method in \cite{weinan1994numerical,guo2015maximum}, the high order SL WENO scheme in \cite{qiu_shu_sl,christlieb2014high,xiong2016conservative,xiong2018high}, the DG method in \cite{liu2000high,zhangshu2010,zhu2017h} and the spectral element method in \cite{fischer2001filter,xu2006stabilization}.
We solve this problem up to $T=8$ and present the surface plots of $\omega$ for both SLDG methods in Figure \ref{figure:shear_contour1}. The solutions by the splitting scheme (left) and non-splitting scheme (middle) are compared against the reference solution (right) computed by the same methods with a doubly refined mesh. In Figure \ref{shear_1d}, we further compare the 1D cuts of the numerical solutions  along $x=\pi$. It is observed that the  non-splitting SLDG method performs much better in correctly resolving solution structures than the splitting counterpart. We also plot the time history of the relative deviations of energy and entropy in Figure \ref{shear_norms}, using a mesh of $100\times100$ cells and $CFL=1$. Again, the non-splitting SLDG method does a better job of conserving the physical invariants than the splitting counterpart.

\begin{figure}[h!]
\centering                              
\includegraphics[height=45mm]{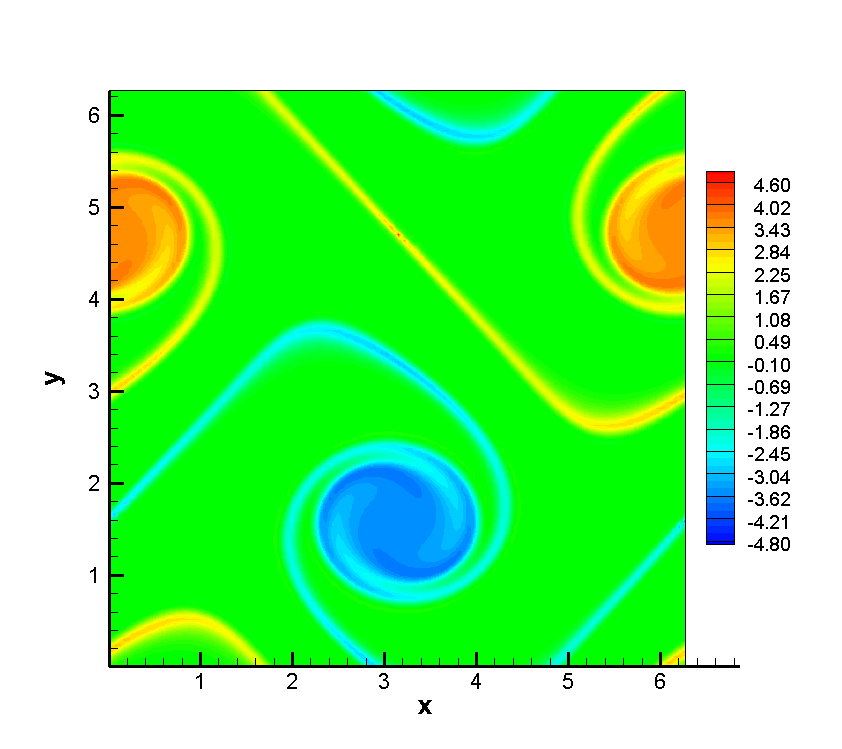}
\includegraphics[height=45mm]{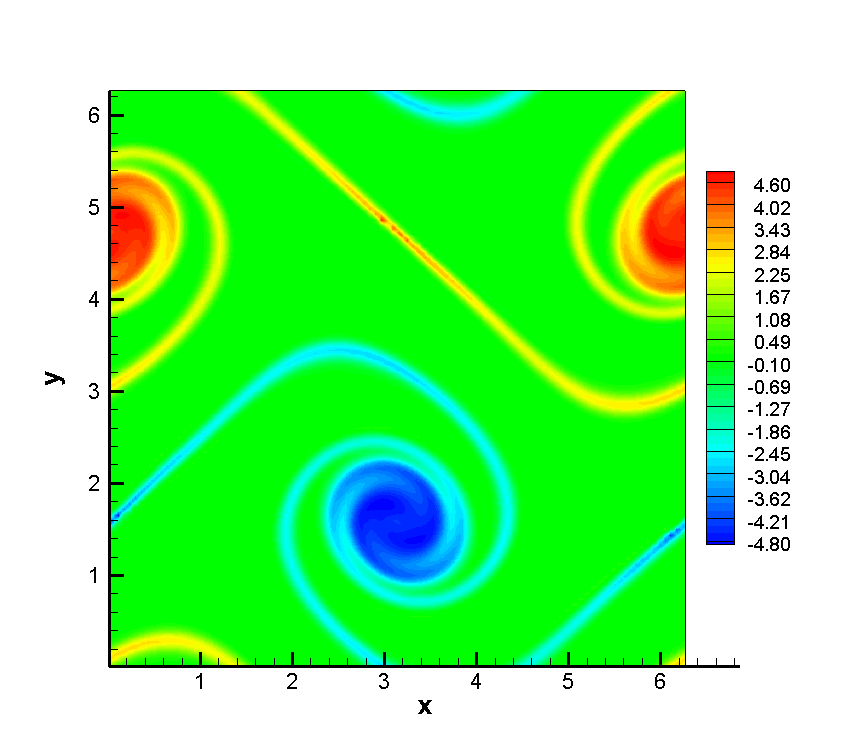}
\includegraphics[height=45mm]{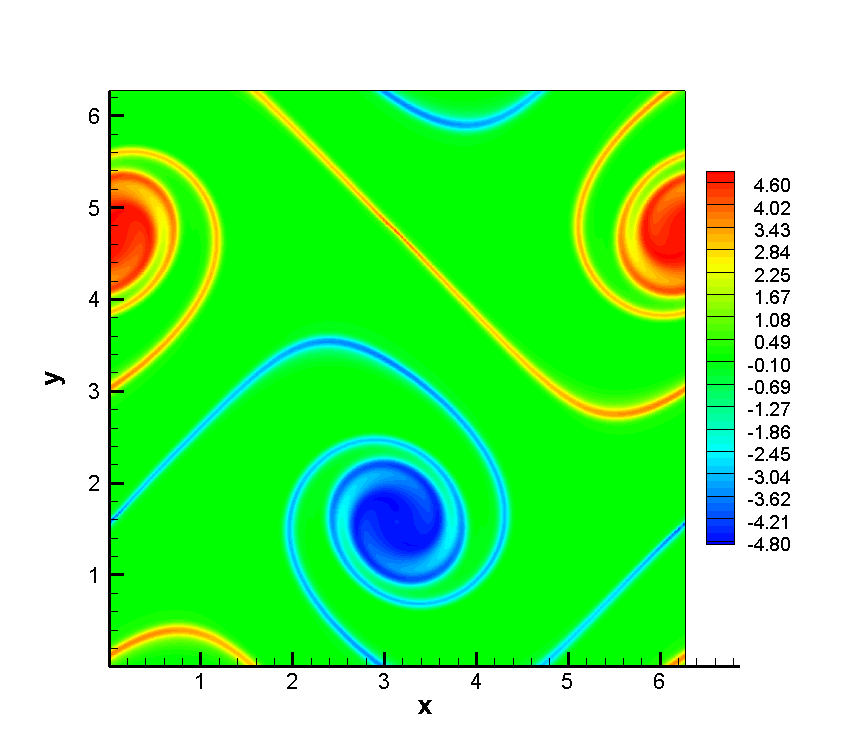}
\includegraphics[height=45mm]{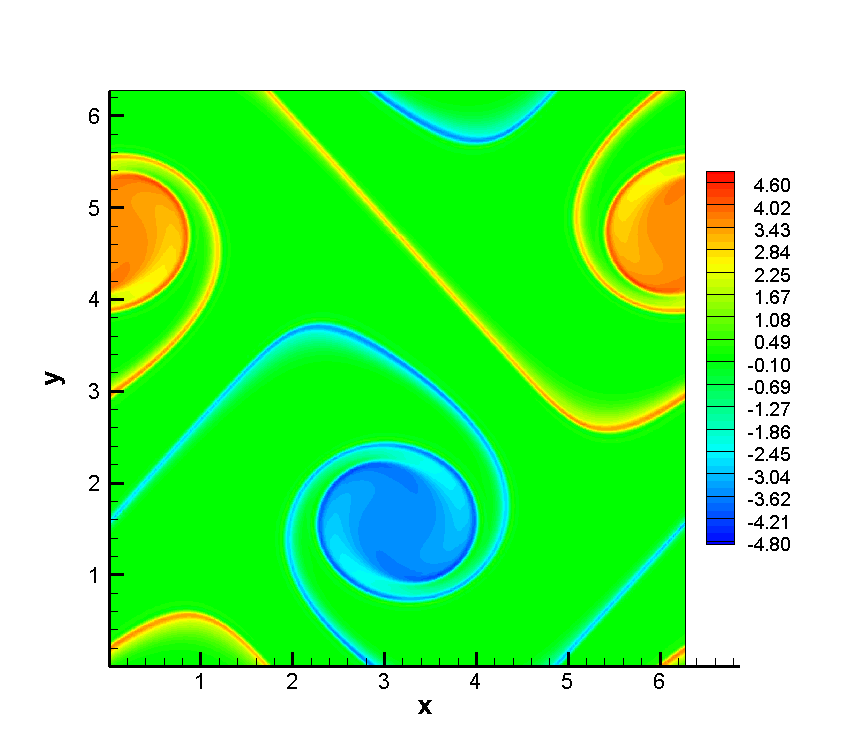}
\includegraphics[height=45mm]{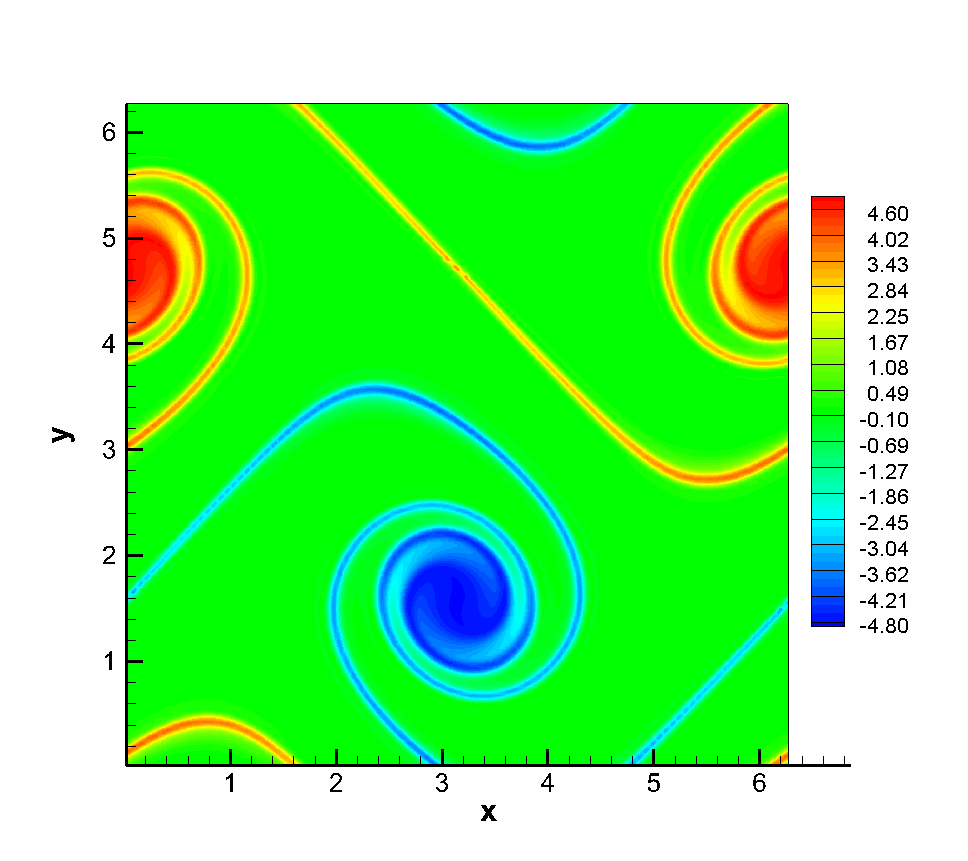}
\includegraphics[height=45mm]{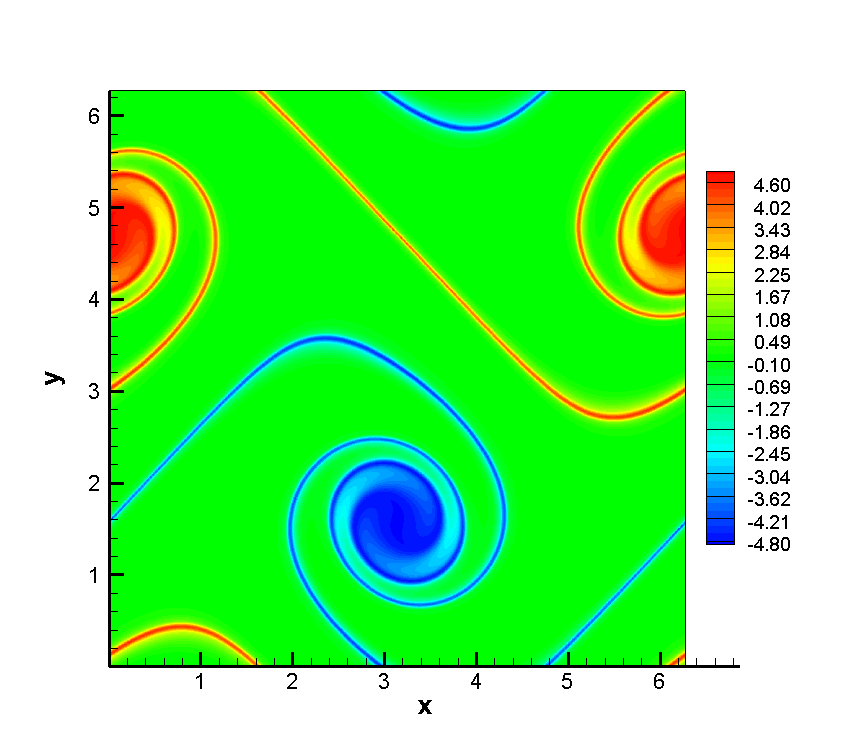}
\caption{Surface plots of the numerical solutions for the shear flow at $T=8$.
A mesh of $100\times100$ cells is used for the left and middle plots, in comparison to the reference solution using a mesh of $200\times200$ cells (right plot).
$CFL=1$.
Upper left: $Q^1$ SLDG-split-$P^1$ LDG.
Upper middle: $P^1$ SLDG-$P^1$ LDG-time2.
Upper right: $P^1$ SLDG-$P^1$ LDG-time2.
Bottom left: $Q^2$ SLDG-split-$P^2$ LDG.
Bottom middle: $P^2$ SLDG-QC-$P^2$ LDG-time3.
Bottom right: $P^2$ SLDG-QC-$P^2$ LDG-time3.
}
\label{figure:shear_contour1}
\end{figure}

\begin{figure}[h!]
\centering                              
\includegraphics[height=70mm]{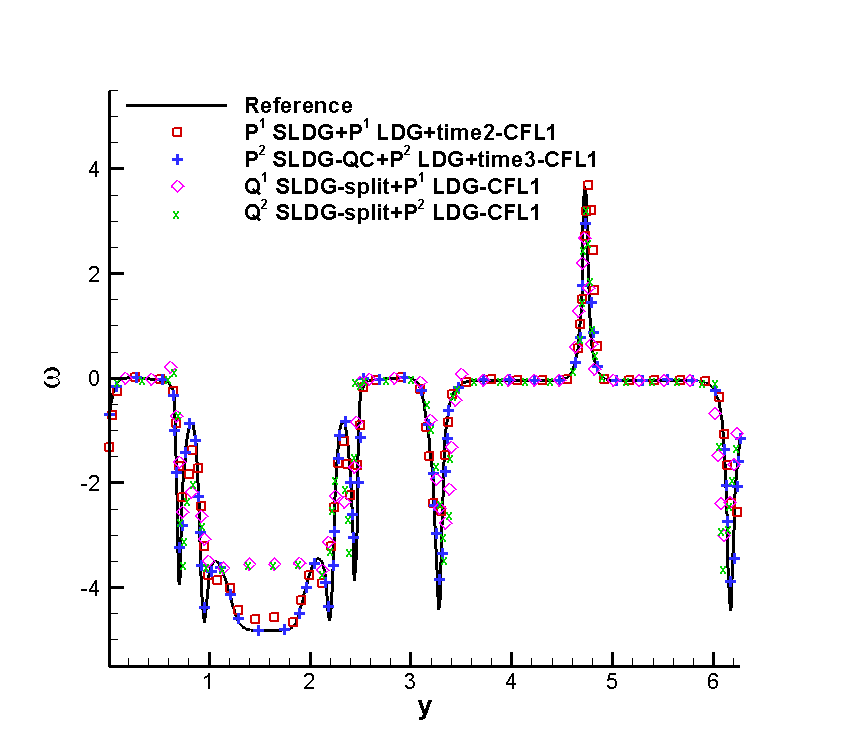}
\caption{1D cut along $x=\pi$ for the shear flow problem at $T=8$. The mesh of $100\times100$ cells is used.
 Lines: A reference solution is obtained by $P^2$ SLDG-QC+$P^2$ LDG+time3 with $CFL=1$ and using a mesh of $200\times200$ cells.
 }
\label{shear_1d}
\end{figure}

\begin{figure}[h!]
\centering                              
\includegraphics[height=70mm]{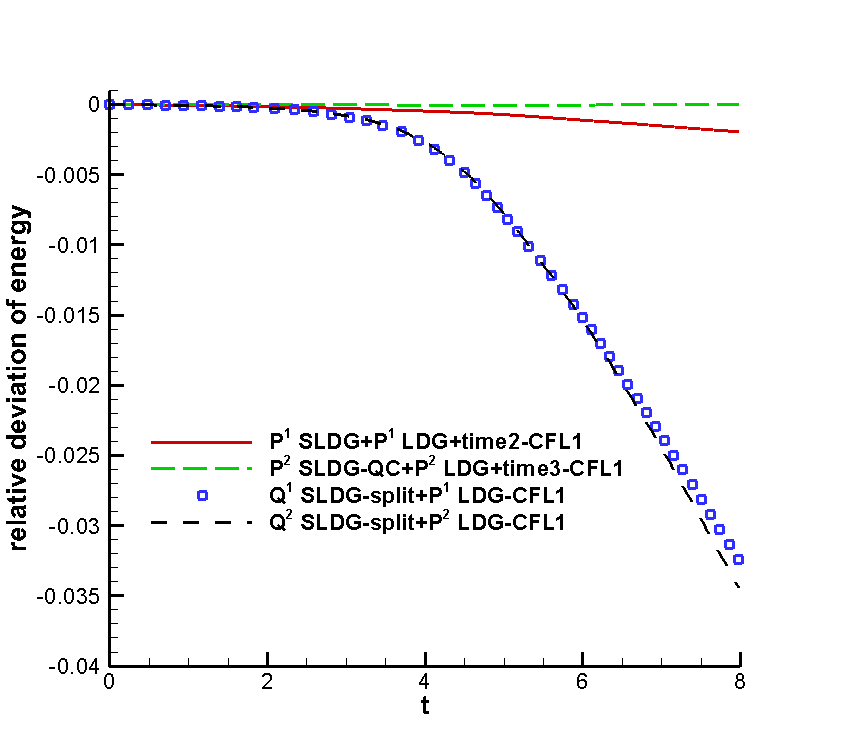}
\includegraphics[height=70mm]{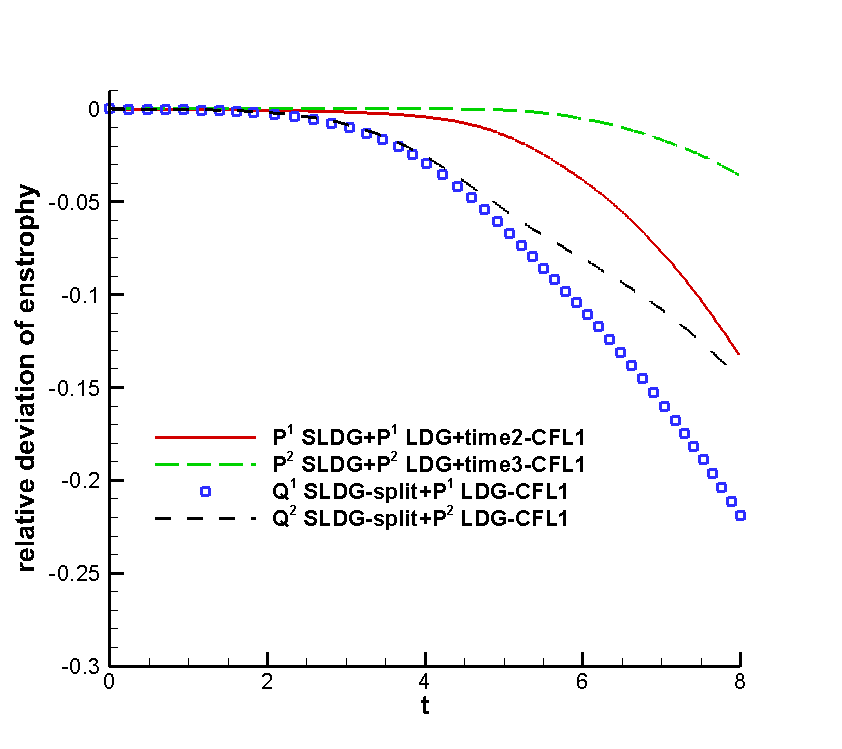}
\caption{Time evolution of the relative deviation of energy (left) and enstrophy (right) for the SLDG method with or without operator splitting for the double shear layer problem. The mesh of $100\times100$ cells is used. }
\label{shear_norms}
\end{figure}

\end{exa}

%


\begin{exa}
(Kelvin-Helmholtz instability.)
The last example is the 2D guiding center model problem with the initial condition
\begin{equation}
\rho_0(x,y) = \sin(y) + 0.015 \cos(k_0x)
\end{equation}
and periodic boundary conditions on the domain $[0,4\pi]\times[0,2\pi]$.
We let $k_0=0.5$, and hence create a Kelvin-Helmholtz instability \cite{shoucri1981two}.
In the literature, this problem was well studied by many authors. 
In Figure \ref{figure:KH_contour1}, we present the surface plots of charge density $\rho$ at $T=40$ by the two SLDG methods with a mesh of $100\times100$ elements. The reference solution with a doubly refined mesh is provided as well.
We further compare the numerical solutions via their 1D cuts along the line $y=\pi$ in Figure \ref{KH_1d}.
As with the previous example, we observe that the solutions by the non-splitting SLDG method  match the reference solution more closely. Furthermore,
the non-splitting method is able to better conserve energy and enstrophy of the system  than the splitting counterpart, see
Figure \ref{KH_norms}.

\begin{figure}[h!]
\centering                              
\includegraphics[height=44mm, width = 50mm]{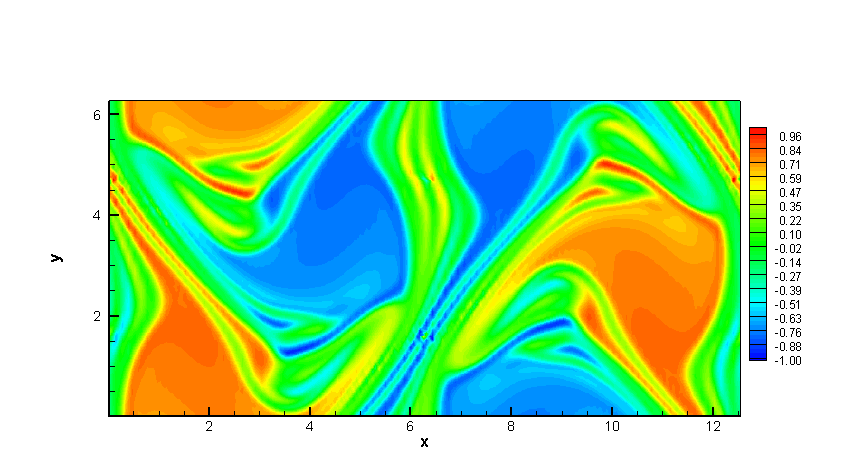}
\includegraphics[height=44mm, width = 50mm]{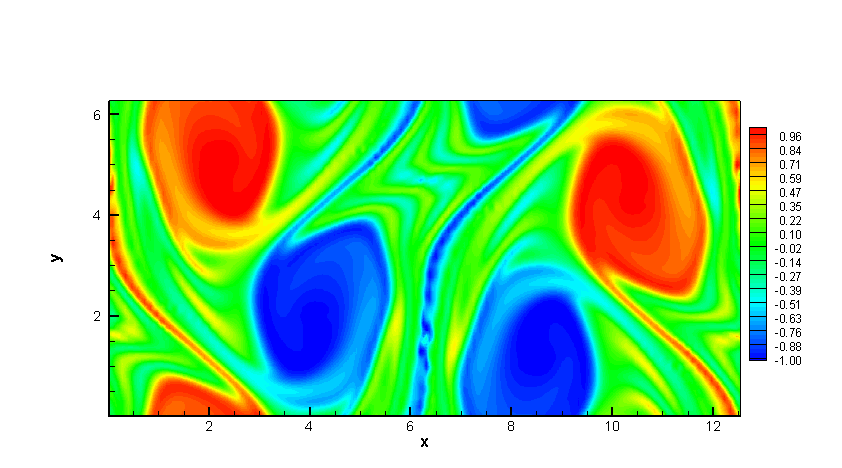}
\includegraphics[height=44mm, width = 50mm]{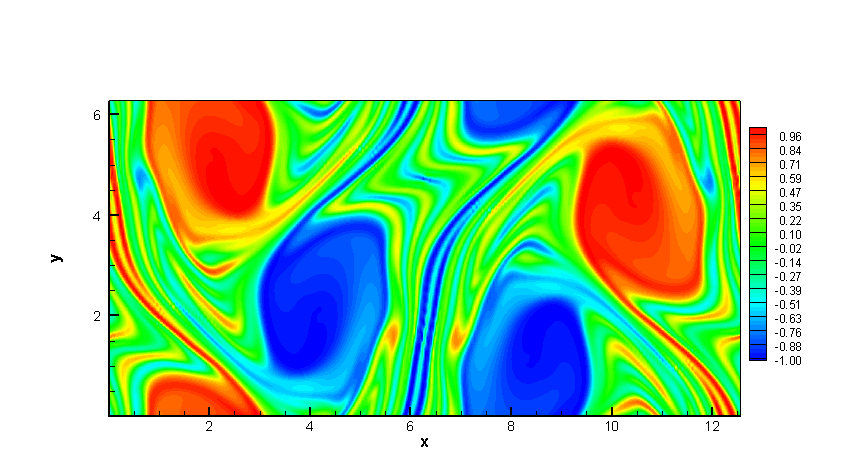}
\includegraphics[height=44mm, width = 50mm]{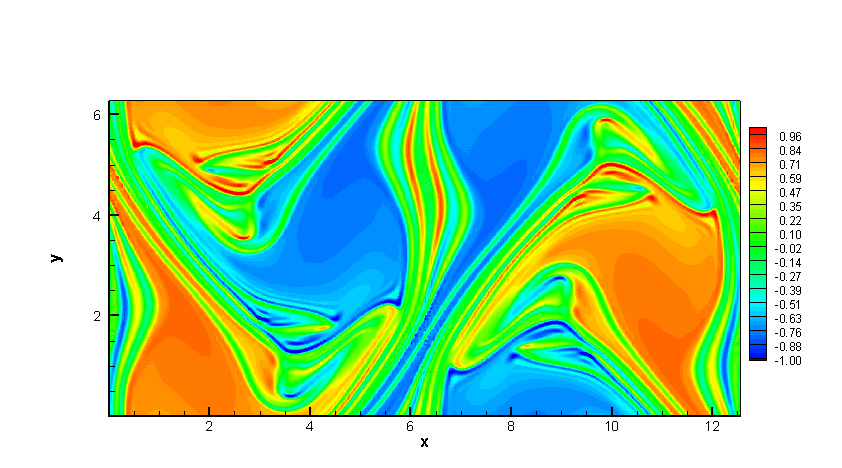}
\includegraphics[height=44mm, width = 50mm]{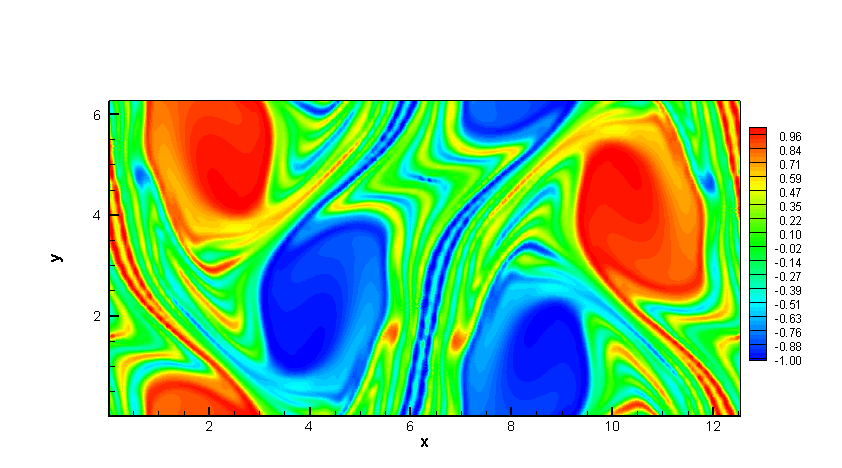}
\includegraphics[height=44mm, width = 50mm]{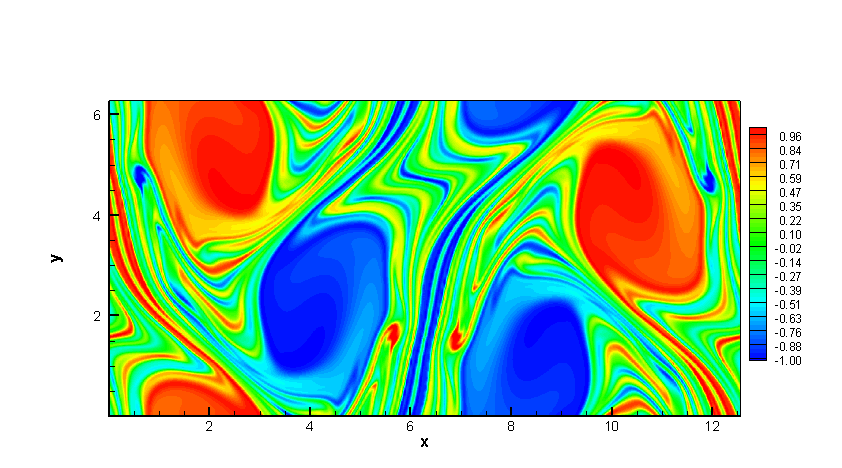}
\caption{Surface plots of the numerical solutions for the Kelvin-Helmholtz instability at $T=40$.
A mesh of $100\times100$ cells is used for the left and middle plots, in comparison to the reference solution using a mesh of $200\times200$ cells (right plot).
$CFL=1$.
Upper left: $Q^1$ SLDG-split-$P^1$ LDG.
Upper middle: $P^1$ SLDG-$P^1$ LDG-time2.
Upper right: $P^1$ SLDG-$P^1$ LDG-time2.
Bottom left: $Q^2$ SLDG-split-$P^2$ LDG.
Bottom middle: $P^2$ SLDG-QC-$P^2$ LDG-time3.
Bottom right: $P^2$ SLDG-QC-$P^2$ LDG-time3.
}
\label{figure:KH_contour1}
\end{figure}

\begin{figure}[h!]
\centering                              
\includegraphics[height=70mm]{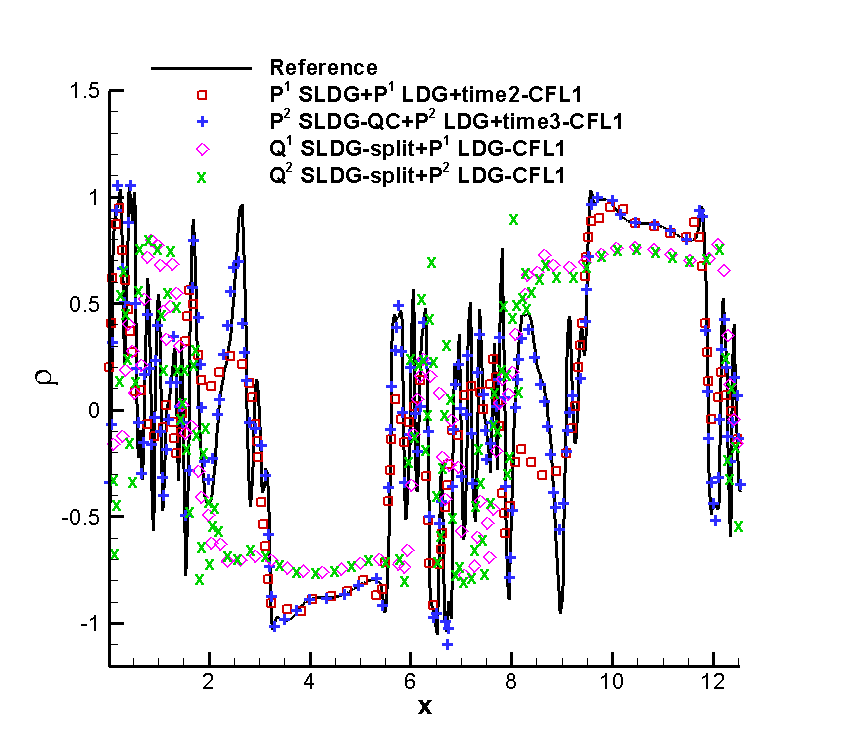}
\caption{1D cut along $y=\pi$ for the Kelvin-Helmholtz instability problem at $T=40$.
 Lines: A reference solution is obtained by $P^2$ SLDG-QC-$P^2$ LDG-time3 with $CFL=1$ and using a mesh of $200\times200$ cells.}
\label{KH_1d}
\end{figure}

\begin{figure}[h!]
\centering                              
\includegraphics[height=70mm]{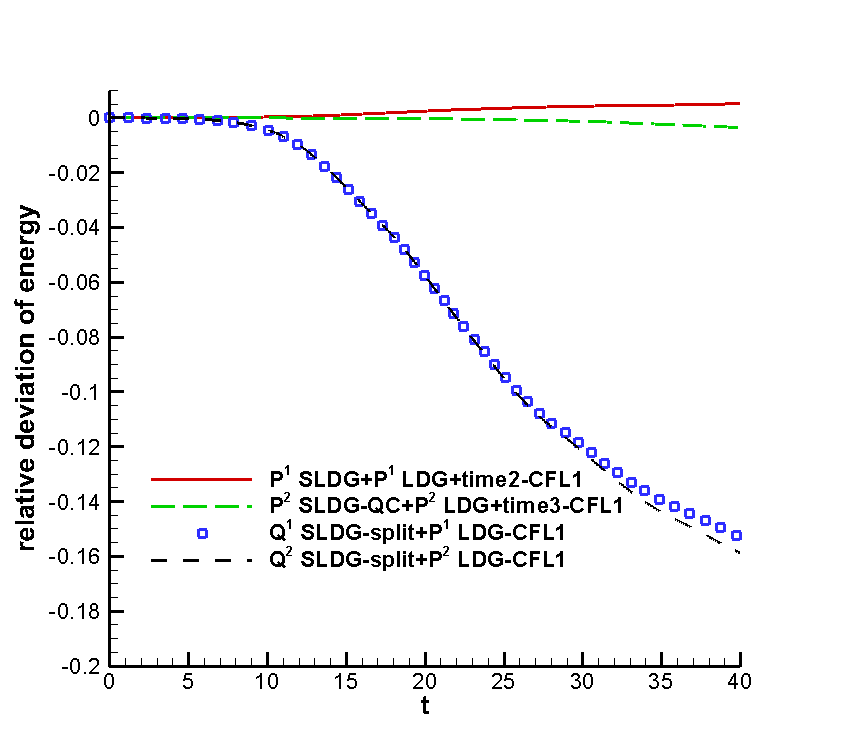}
\includegraphics[height=70mm]{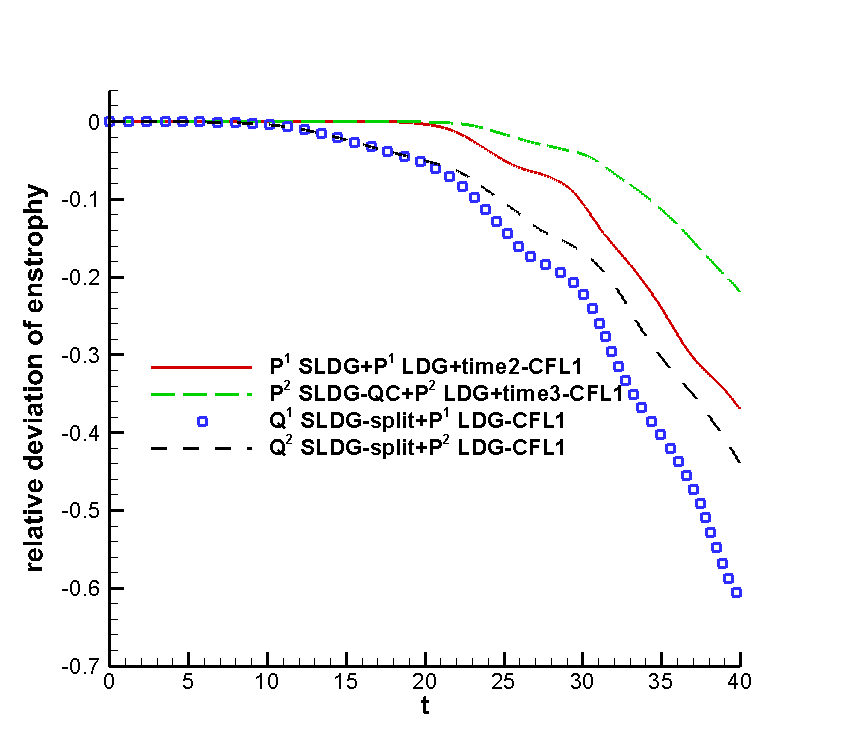}
\caption{Time evolution of the relative deviation of enstrophy for the SLDG method with or without operator splitting for the  Kelvin-Helmholtz instability problem. The mesh $100\times100$ cells of is used.}
\label{KH_norms}
\end{figure}

\end{exa}

\section{Conclusion}
\label{sec4}
\setcounter{equation}{0}
\setcounter{figure}{0}
\setcounter{table}{0}

In this paper, we review several SLDG methods and compare their performance for solving 2D linear transport problems, and nonlinear Vlasov and incompressible Euler equations. The performances of SLDG  schemes in terms of error magnitude, order of accuracy, CPU cost, resolution of complicated solution structures, as well as their ability to conserve important physical invariants were benchmarked through extensive numerical experiments. In general, when the geometry is simple, and smaller $CFL$ is needed for accuracy, the scheme based on dimensional splitting together with 1D SLDG formulation is preferred for its simplicity in implementation; when larger $CFL$ is desired for efficiency without significantly sacrificing accuracy, the truly multi-dimensional SLDG scheme is preferred for its efficiency and efficacy in resolving solutions using a relatively coarse mesh and extra large time stepping sizes.


\bibliographystyle{abbrv}
\bibliography{refer17}

\end{document}